\documentclass[11pt]{article}

\usepackage{amssymb,latexsym,amsmath}

\usepackage{graphicx}

\begin{document}

\newcommand{\bfi}{\bfseries\itshape}

\makeatletter

\@addtoreset{figure}{section}

\def\thefigure{\thesection.\@arabic\c@figure}

\def\fps@figure{h, t}

\@addtoreset{table}{bsection}

\def\thetable{\thesection.\@arabic\c@table}

\def\fps@table{h, t}

\@addtoreset{equation}{section}

\def\theequation{\thesubsection.\arabic{equation}}

\makeatother

\newtheorem{thm}{Theorem}[section]

\newtheorem{prop}[thm]{Proposition}

\newtheorem{lema}[thm]{Lemma}

\newtheorem{cor}[thm]{Corollary}

\newtheorem{defi}[thm]{Definition}

\newtheorem{rk}[thm]{Remark}

\newtheorem{exempl}{Example}[section]

\newenvironment{exemplu}{\begin{exempl}  \em}{\hfill $\square$

\end{exempl}}

\newcommand{\comment}[1]{\par\noindent{\raggedright\texttt{#1}

\par\marginpar{\textsc{Comment}}}}

\newcommand{\todo}[1]{\vspace{5 mm}\par \noindent \marginpar{\textsc{ToDo}}\framebox{\begin{minipage}[c]{0.95 \textwidth}

\tt #1 \end{minipage}}\vspace{5 mm}\par}

\newcommand{\ea}{\mbox{{\bf a}}}

\newcommand{\eu}{\mbox{{\bf u}}}

\newcommand{\ueu}{\underline{\eu}}

\newcommand{\ueo}{\overline{u}}

\newcommand{\oeu}{\overline{\eu}}

\newcommand{\ew}{\mbox{{\bf w}}}

\newcommand{\ef}{\mbox{{\bf f}}}

\newcommand{\eF}{\mbox{{\bf F}}}

\newcommand{\eC}{\mbox{{\bf C}}}

\newcommand{\en}{\mbox{{\bf n}}}

\newcommand{\eT}{\mbox{{\bf T}}}

\newcommand{\eL}{\mbox{{\bf L}}}

\newcommand{\eR}{\mbox{{\bf R}}}

\newcommand{\eV}{\mbox{{\bf V}}}

\newcommand{\eU}{\mbox{{\bf U}}}

\newcommand{\ev}{\mbox{{\bf v}}}

\newcommand{\eve}{\mbox{{\bf e}}}

\newcommand{\uev}{\underline{\ev}}

\newcommand{\eY}{\mbox{{\bf Y}}}

\newcommand{\eK}{\mbox{{\bf K}}}

\newcommand{\eP}{\mbox{{\bf P}}}

\newcommand{\eS}{\mbox{{\bf S}}}

\newcommand{\eJ}{\mbox{{\bf J}}}

\newcommand{\eB}{\mbox{{\bf B}}}

\newcommand{\eH}{\mbox{{\bf H}}}

\newcommand{\leb}{\mathcal{ L}^{n}}

\newcommand{\eI}{\mathcal{ I}}

\newcommand{\eE}{\mathcal{ E}}

\newcommand{\hen}{\mathcal{H}^{n-1}}

\newcommand{\eBV}{\mbox{{\bf BV}}}

\newcommand{\eA}{\mbox{{\bf A}}}

\newcommand{\eSBV}{\mbox{{\bf SBV}}}

\newcommand{\eBD}{\mbox{{\bf BD}}}

\newcommand{\eSBD}{\mbox{{\bf SBD}}}

\newcommand{\ecs}{\mbox{{\bf X}}}

\newcommand{\eg}{\mbox{{\bf g}}}

\newcommand{\paromega}{\partial \Omega}

\newcommand{\gau}{\Gamma_{u}}

\newcommand{\gaf}{\Gamma_{f}}

\newcommand{\sig}{{\bf \sigma}}

\newcommand{\gac}{\Gamma_{\mbox{{\bf c}}}}

\newcommand{\deu}{\dot{\eu}}

\newcommand{\dueu}{\underline{\deu}}

\newcommand{\dev}{\dot{\ev}}

\newcommand{\duev}{\underline{\dev}}

\newcommand{\weak}{\stackrel{w}{\approx}}

\newcommand{\mild}{\stackrel{m}{\approx}}

\newcommand{\strong}{\stackrel{s}{\approx}}

\newcommand{\weakdown}{\rightharpoondown}

\newcommand{\opg}{\stackrel{\mathfrak{g}}{\cdot}}

\newcommand{\opunu}{\stackrel{1}{\cdot}}
\newcommand{\opdoi}{\stackrel{2}{\cdot}}

\newcommand{\opn}{\stackrel{\mathfrak{n}}{\cdot}}
\newcommand{\opx}{\stackrel{x}{\cdot}}

\newcommand{\tr}{\ \mbox{tr}}

\newcommand{\Ad}{\ \mbox{Ad}}

\newcommand{\ad}{\ \mbox{ad}}

\renewcommand{\contentsname}{ }

\title{Dilatation structures II. \\  Linearity, self-similarity and the 
Cantor set}

\author{Marius Buliga \\
\\
Institute of Mathematics, Romanian Academy \\
P.O. BOX 1-764, RO 014700\\
Bucure\c sti, Romania\\
{\footnotesize Marius.Buliga@imar.ro}}

\date{This version:  14.02.2007}

\maketitle

\section*{Introduction}

In this paper we continue the study of dilatation structures, introduced in  \cite{buligadil1}. 

A dilatation structure on a metric space is a kind of enhanced self-similarity. 
By way of examples this is explained here with the help of the middle-thirds Cantor set.

Linear and self-similar dilatation structures are introduced  and 
studied on ultrametric spaces, especially on the boundary of the dyadic tree 
(same as the middle-thirds Cantor set). 

Some other  examples of dilatation structures, which share some common 
features, are given. Another class of examples, coming from sub-Riemannian 
geometry, will make the subject of an article in preparation. 

In the particular case of ultrametric spaces the axioms of dilatation structures 
take a simplified form, leading to a description of all possible weak dilatation 
structures on the Cantor set.

 As an  application we prove that there is more than  one linear and self-similar 
dilatation structure on the Cantor set, compatible with  the iterated 
functions system which defines the Cantor set.

Applications to self-similar groups are reserved for a further paper.

\tableofcontents

\section{Basics}

Here we collect well-known facts and notations which we shall need further. 

\subsection{Notations}

Let $\Gamma$ be  a topological separated commutative group  endowed with a continuous group morphism 
$$\nu : \Gamma \rightarrow (0,+\infty)$$ with $\displaystyle \inf \nu(\Gamma)  =  0$. Here $(0,+\infty)$ is 
taken as a group with multiplication. The neutral element of $\Gamma$ is denoted by $1$. We use the multiplicative notation for the operation in $\Gamma$. 

The morphism $\nu$ defines an invariant topological filter on $\Gamma$ (equivalently, an end). Indeed, 
this is the filter generated by the open sets $\displaystyle \nu^{-1}(0,a)$, $a>0$. From now on 
we shall name this topological filter (end) by "0" and we shall write $\varepsilon \in \Gamma \rightarrow 
0$ for $\nu(\varepsilon)\in (0,+\infty) \rightarrow 0$. 

The set $\displaystyle \Gamma_{1} = \nu^{-1}(0,1] $ is a semigroup. We note $\displaystyle 
\bar{\Gamma}_{1}= \Gamma_{1} \cup \left\{ 0\right\}$
On the set $\displaystyle 
\bar{\Gamma}= \Gamma \cup \left\{ 0\right\}$ we extend the operation on $\Gamma$ by adding the rules  
$00=0$ and $\varepsilon 0 = 0$ for any $\varepsilon \in \Gamma$. This is in agreement with the invariance 
of the end $0$ with respect to translations in $\Gamma$.

\subsection{Words and the Cantor middle-thirds set}

Let $X$ be a finite, non empty  set. The elements of $X$ are called letters.  
The collection of words of finite length in the alphabet $X$ is denoted by $X^{*}$. The empty word 
$\emptyset$ is an element of $X^{*}$. 

In this paper we shall work mainly with the alphabet 
$\displaystyle X = \left\{ 0, 1 \right\}$. In this particular case one can define a 
conjugate function (from $X$ to $x$), given by $\displaystyle \bar{0} = 1$, 
$\displaystyle \bar{1} = 0$.

 The length of any word $w\in X^{*}$ 
$$w = a_{1} ... a_{m}  \, , \quad a_{k} \in X  \quad \forall k=1, ... , m $$
is denoted by $\mid w \mid = m$. 

The set of words which are infinite at right is denoted by 
$$X^{\omega} = \left\{ f \  \mid \ \   f: \mathbb{N}^{*}\rightarrow X \right\} = X^{\mathbb{N}^{*}} \quad . $$
Concatenation of words is naturally defined. If $\displaystyle q_{1},q_{2} \in X^{*}$ and 
$w \in X^{\omega}$ then $\displaystyle q_{1}q_{2} \in X^{*}$ and $q_{1}w \in X^{\omega}$. 

The shift map $\displaystyle s : X^{\omega} \rightarrow X^{\omega}$ is defined by 
$$ w = w_{1} \,Ês(w) $$
for any word $\displaystyle w \in X^{\omega}$. For any $\displaystyle k \in \mathbb{N}^{*}$ we define  
$\displaystyle [w]_{k} \in X^{k} \subset X^{*}$, $\displaystyle \left\{ w \right\}_{k} \in X^{\omega}$ by 
$$ w =  [w]_{k} \,Ês^{k}(w)  \quad , \quad \left\{ w \right\}_{k}  = s^{k}(w) \quad . $$

The topology on $\displaystyle X^{\omega}$ is generated by cylindrical sets $qX^{\omega}$, 
for all $q\in X^{*}$.  The topological space $\displaystyle X^{\omega}$  is compact . 

To any $q\in X^{*}$ is associated a continuous injective transformation $\hat{q}:X^{\omega}\rightarrow 
X^{\omega}$, $\hat{q}(w) = qw$. The semigroup $X^{*}$ (with respect to concatenation) can be identified with the semigroup (with respect to function composition) of these transformations. This semigroup is obviously generated by $X$. The empty word corresponds to the identity function. 

The topological space $X^{\omega}$ is metrizable. Indeed, denote by $N$ the number of letters 
in the alphabet $X$. Let us fix a bijection between $X$ and $\left\{1, ... , N\right\}$, which allows us to 
identify $X$ with $\left\{1, ... , N\right\}$. 

For any $w\in X^{\omega}$ and $i\in\mathbb{N}^{*}$, the symbol $\displaystyle w_{i}\in X=\left\{1, ... , N\right\}$ denotes 
the $i$-th letter in the word $w$. Define then  the injective continuous function 
$$\Phi: X^{\omega}\rightarrow \mathbb{R} \quad , \quad \Phi(w) = \sum_{i=1}^{\infty} \frac{w_{i}}{2^{i}} \quad . $$

If $N$ is a prime number then we may equally define a homeomorphism between 
$X^{\omega}$ and the group of $N$-adic integers. We shall first recall some basic facts 
about p-adic integers.  

$\displaystyle \mathbb{Q}_{p}$ is the closure of $\displaystyle \mathbb{Q}$ with respect to 
$p$-adic norm. For $x\in \mathbb{Z}$, $x\not = 0$ the $p$-norm is defined by 
$$\mid x \mid_{p} = \inf \left\{ p^{-r} \, \mid \, p^{r} \mbox{ divides } x \right\} \quad . $$
For $x/y\in \mathbb{Q}$ ($x,y\in \mathbb{Z}$,  $y\not = 0$) the $p$-adic norm is defined by 
$$\mid \frac{x}{y}\mid_{p} = \frac{\mid x \mid_{p}}{\mid y \mid_{p}} \quad .$$
The norm induces an ultrametric  distance on $\displaystyle \mathbb{Q}_{p}$. The closed unit disk of $\displaystyle \mathbb{Q}_{p}$ is called the ring of $p$-adic integers, denoted $\displaystyle \mathbb{Z}_{p}$. The set  $\displaystyle \mathbb{Z}_{p}$ is compact. 

\begin{prop}
Any element of  $\displaystyle x\in \mathbb{Z}_{p}$ admits an unique $p$-adic expansion 
$$x = \sum_{i=1}^{\infty} x_{i} p^{i} \quad , $$
with all $x_{i}\in \left\{0, ... , p-1\right\}$.  
Any element of  $\displaystyle x\in \mathbb{Q}_{p}$ admits an unique $p$-adic expansion 
$$x = \sum_{i=r}^{\infty} x_{i} p^{i} \quad , $$
starting from some $r\in \mathbb{Z}$, with all $x_{i}\in \left\{0, ... , p-1\right\}$.  
\end{prop}

The addition and multiplication of $p$adic numbers is done using the $p$-adic expansion and the 
standard algorithms (with some simple modifications: with addition or multiplications of "digits" $\displaystyle x_{i}$, modulo $p$, with remainders, from left to right).

The function 
$$\Psi: X^{\omega}\rightarrow \mathbb{Z}_{N} \quad , \quad \Phi(w) = \sum_{i=1}^{\infty} (w_{i}-1)N^{i}  $$
is then  a homeomorphism. 

The elements of $X^{*}$, seen as transformations of $X^{\omega}$, are contractions, with respect to the distances induced by $\Phi$ and ( if $N$ is prime) $\Psi$. Indeed, it is enough to check this for the transformations associated with the letters in the alphabet $X$. Let $a\in X$ and $w\in X^{\omega}$. 
We have then 
$$ \Phi(aw) = \frac{a}{2} + \frac{1}{2}\Phi(w) \quad , \quad \Psi(aw) = (a-1) + N \Psi(w) \quad , $$
which implies that for any $\displaystyle w_{1}, w_{2}\in X^{\omega}$ we have 
$$\mid \Phi(aw_{1}) - \Phi(aw_{2}) \mid = \frac{1}{2} \mid \Phi(w_{1}) -\Phi(w_{2})\mid \quad , $$ 
$$\mid \Psi(aw_{1}) - \Psi(aw_{2}) \mid_{N} \,                                               =  \,Ê                                        \frac{1}{N} \mid \Phi(w_{1}) -\Phi(w_{2})\mid_{N} \quad . $$
Therefore $\displaystyle \Phi(X^{\omega}), \Psi(X^{\omega})$ are invariant sets of iterated functions systems of contractions. 

In the particular case $\displaystyle X = \left\{ 0, 1 \right\}$, up to a 
multiplicative factor $\Phi(X^{\omega})$ is the middle-thirds  Cantor set.

\subsection{IFS of contractions}

In order to put things into perspective, we shall recall simple facts about iterated functions systems of contractions and their invariant sets, following Hutchinson \cite{hutch}. 

\begin{defi}
A contraction is a Lipschitz map $\phi:(X,d)\rightarrow (X,d)$ with Lipschitz constant smaller than $1$. 

An iterated system of contractions $\mathcal{S}$ is a finite collection of contractions on a complete 
metric space $(X,d)$. 

An invariant set of $\mathcal{S}$ is a set $M\subset X$ such that  
$$M = \bigcup_{\phi\in \mathcal{S}} \phi(M) \quad .$$
\end{defi}

\begin{thm}
There exists and it is unique  a non empty  bounded invariant set of $\mathcal{S}$, denoted by $K(\mathcal{S})$. Moreover $K(\mathcal{S})$ is compact. 

For any bounded, non empty set $A \subset X$, let us define 
$$\mathcal{S}(A) = \bigcup_{\phi\in \mathcal{S}} \phi(A) \quad .$$ 
Then $\displaystyle \mathcal{S}^{n}(A)$ converges in the Hausdorff distance to $K(\mathcal{S})$, as 
$n\rightarrow \infty$. 
\label{tifs}
\end{thm}

If $\phi:X\rightarrow X$ is a contraction and $(X,d)$ is compact then it has an unique fixed point 
$\displaystyle x_{\phi}\in X$, that is $\displaystyle x_{\phi}$  exists and it is unique with the property 
$\displaystyle \phi(x_{\phi}) =  x_{\phi}$.  

Let $\mathcal{S}^{*}$ be the semigroup (with function composition) generated by $\mathcal{S}$ and 
$\displaystyle Fix(\mathcal{S})$ be the collection of fixed points of elements of $\mathcal{S}^{*}$ (recall 
that each element of  $\mathcal{S}^{*}$ is a contraction which preserves the compact  set $M(\mathcal{S})$). 

\begin{thm}
The set $\displaystyle Fix(\mathcal{S})$ is dense in $M(\mathcal{S})$.
\label{thutch} 
\end{thm}

We can give codes to elements of $M(\mathcal{S})$. Indeed, let us start by 
remarking that we already used a notation similar to one in the previous 
subsection, namely $\mathcal{S}^{*}$. There 
is a surjective morphism from $\mathcal{S}^{*}$, as the semigroup of finite 
words with concatenation, to 
$\mathcal{S}^{*}$, as the semigroup of contractions generated by $\mathcal{S}$, 
with function composition.  This morphism induces a surjective function 
$\displaystyle 
\Lambda: \mathcal{S}^{\omega} \rightarrow K(\mathcal{S})$ with the property 
that for any 
$\displaystyle q=\phi_{1} ... \phi_{m} \in \mathcal{S}^{*}$ (finite word) 
and any 
$\displaystyle w\in \mathcal{S}^{\omega}$ we have 
$$\Lambda(qw) = \phi_{1} ... \phi_{m} ( \Lambda(w)) \quad . $$

The function $\Lambda$ is constructed like this:  for any 
$\displaystyle w\in \mathcal{S}^{\omega}$ and $n\in \mathbb{N}^{*}$, let 
$\displaystyle [w]_{n} \in \mathcal{S}^{*}$ be the $n$-letter word from the 
beginning of $w$. There exists $\displaystyle w' \in \mathcal{S}^{\omega}$ 
such that $\displaystyle w = [w]_{n} w'$.  Define $\Lambda(w)$ by 
$$\Lambda(w) \in \bigcap_{n \in \mathbb{N}^{*}} [w]_{n}(K(\mathcal{S})) 
\quad , $$
(where $\displaystyle [w]_{n}$ from the right hand side of the previous 
relation is understood as the composition of the first $n$ letters of the 
word $w$). 

The definition is good because $\displaystyle diam \, [w]_{n}(K(\mathcal{S})) 
\rightarrow 0$ as $n\rightarrow \infty$ and for any 
$\displaystyle n\in \mathbb{N}^{*}$ we have 
$$[w]_{n+1}(K(\mathcal{S})) \subset [w]_{n}(K(\mathcal{S})) \quad .$$ 
Therefore the intersection of all 
$\displaystyle [w]_{n}(K(\mathcal{S}))$ is a singleton. 

The function $\Lambda$ is not generally bijective. It is, though, if the following condition is satisfied. 

\begin{defi}
$\mathcal{S}$ satisfies the open set condition if there exists a non empty open set $U$ such that 
\begin{enumerate}
\item[(a)] $\displaystyle \bigcup_{\phi\in \mathcal{S}} \phi(U) \subset U$, 
\item[(b)] if $\phi\not = \psi$, $\phi, \psi\in \mathcal{S}$, then $\phi(U)\cap\psi(U) = \emptyset$. 
\end{enumerate}
\label{defopen}
\end{defi}

\subsection{Isometries of  the dyadic tree}

The dyadic tree $\mathcal{T}$  is the infinite rooted binary tree, 
with any node having two descendants. The nodes are  coded by elements of 
$\displaystyle X^{*}$, $X = \left\{ 0,1\right\}$. The root is coded by the 
empty word $\emptyset$  and if a node is coded by $x\in X^{*}$ then its left 
hand side descendant has the code $x0$ and its 
right hand side descendant has the code $x1$. We shall therefore identify the 
dyadic tree with $\displaystyle X^{*}$ and we put on the dyadic tree the 
natural (ultrametric) distance on $\displaystyle X^{*}$. The boundary 
(or the set of ends) of the dyadic tree is then the same as the compact 
ultrametric space $\displaystyle X^{\omega}$. 

An isometry of $\mathcal{T}$ is just an invertible transformation which 
preserves the structure of the 
tree. It is well known that isometries  of $\displaystyle (X^{\omega}, d)$ are 
the same as isometries of $\mathcal{T}$.  

Let $\displaystyle A \in Isom(X^{\omega}, d)$ be such an isometry. For any 
finite word $\displaystyle q\in X^{*}$ we may define $\displaystyle A_{q} \in 
Isom(X^{\omega}, d)$ by $$A(qw) = A(q) \, A_{q}(w)$$
for any $\displaystyle w \in X^{\omega}$. Note that in the previous relation 
$A(q)$ makes sense because $A$ is also an isometry of $\mathcal{T}$.

\section{Dilatation structures}

The first two  sections contain  notions and results introduced or proved in 
Buliga \cite{buligadil1}. 
The space $(X,d)$ is a complete, locally compact metric space.

\subsection{Axioms of dilatation structures}
The  axioms of  a dilatation structure $(X,d,\delta)$ are listed further. 
The first axiom is merely a preparation for the next axioms. That is why we 
counted it as axiom 0.

\begin{enumerate}
\item[{\bf A0.}] The dilatations $$ \delta_{\varepsilon}^{x}: U(x) 
\rightarrow V_{\varepsilon}(x)$$ are defined for any 
$\displaystyle \varepsilon \in \Gamma, \nu(\varepsilon)\leq 1$. 
All dilatations are homeomorphisms (invertible, continuous, with 
continuous inverse). 

We suppose  that there is  $1<A$ such that for any $x \in X$ we have 
$$\bar{B}_{d}(x,A) \subset U(x)  \ .$$
 We suppose that for all $\varepsilon \in \Gamma$, $\nu(\varepsilon) \in 
(0,1)$, we have 
$$ B_{d}(x,\varepsilon) \subset \delta_{\varepsilon}^{x} B_{d}(x,A) \subset V_{\varepsilon}(x) \subset U(x) \ .$$

 For $\nu(\varepsilon) \in (1,+\infty)$ the associated dilatation  
$$\delta^{x}_{\varepsilon} : W_{\varepsilon}(x) \rightarrow B_{d}(x,B) \ , $$
is injective, invertible on the image. We shall suppose that $\displaystyle  W_{\varepsilon}(x)$ is open, 
$$V_{\varepsilon^{-1}}(x) \subset W_{\varepsilon}(x)$$
and that for all $\displaystyle \varepsilon \in \Gamma_{1}$ and $\displaystyle u \in U(x)$ we have
$$\delta_{\varepsilon^{-1}}^{x} \ \delta^{x}_{\varepsilon} u \ = \ u \ .$$
\end{enumerate}

We remark that we have the following string of inclusions, for any $\varepsilon \in \Gamma$, $\nu(\varepsilon) \leq 1$, and any $x \in X$:
$$ B_{d}(x,\varepsilon) \subset \delta^{x}_{\varepsilon}  B_{d}(x, A) \subset V_{\varepsilon}(x) \subset 
W_{\varepsilon^{-1}}(x) \subset \delta_{\varepsilon}^{x}  B_{d}(x, B) \quad . $$

A further technical condition on the sets  $\displaystyle V_{\varepsilon}(x)$ and $\displaystyle W_{\varepsilon}(x)$  will be given just before the axiom A4. (This condition will be counted as part of 
axiom A0.)

\begin{enumerate}
\item[{\bf A1.}]  We  have 
$\displaystyle  \delta^{x}_{\varepsilon} x = x $ for any point $x$. We also have $\displaystyle \delta^{x}_{1} = id$ for any $x \in X$.

Let us define the topological space
$$ dom \, \delta = \left\{ (\varepsilon, x, y) \in \Gamma \times X \times X \mbox{ : } \quad \mbox{ if } \nu(\varepsilon) \leq 1 \mbox{ then } y \in U(x) \mbox{ , else } y \in W_{\varepsilon}(x) \right\} \quad , $$ 
with the topology inherited from the product topology on $\Gamma \times X \times X$. Consider also 
$\displaystyle Cl(dom \, \delta)$, the closure of $dom \, \delta$ in $\displaystyle \bar{\Gamma} \times X \times X$ with product topology. 
The function 
$$\delta : dom \, \delta \rightarrow  X$$ defined by 
$\displaystyle \delta (\varepsilon,  x, y)  = \delta^{x}_{\varepsilon} y$ is continuous. Moreover, it can be continuously extended to $\displaystyle Cl(dom \, \delta)$ and we have 
$$\lim_{\varepsilon\rightarrow 0} \delta_{\varepsilon}^{x} y \, = \, x \quad . $$

\item[{\bf A2.}] For any  $x, \in K$, $\displaystyle \varepsilon, \mu \in \Gamma_{1}$ and $\displaystyle u \in 
\bar{B}_{d}(x,A)$   we have: 
$$ \delta_{\varepsilon}^{x} \delta_{\mu}^{x} u  = \delta_{\varepsilon \mu}^{x} u  \ .$$

\item[{\bf A3.}]  For any $x$ there is a  function $\displaystyle (u,v) \mapsto d^{x}(u,v)$, defined for any $u,v$ in the closed ball (in distance d) $\displaystyle 
\bar{B}(x,A)$, such that 
$$\lim_{\varepsilon \rightarrow 0} \quad \sup  \left\{  \mid \frac{1}{\varepsilon} d(\delta^{x}_{\varepsilon} u, \delta^{x}_{\varepsilon} v) \ - \ d^{x}(u,v) \mid \mbox{ :  } u,v \in \bar{B}_{d}(x,A)\right\} \ =  \ 0$$
uniformly with respect to $x$ in compact set. 

\end{enumerate}

\begin{rk}
The "distance" $d^{x}$ can be degenerated. That means: there might be 
$\displaystyle v,w \in \bar{B}_{d}(x,A)$ such that $\displaystyle d^{x}(v,w) = 0$ but $v \not = w$. We shall 
use further the name "distance" for $d^{x}$, essentially by commodity, but  keep 
in mind the possible degeneracy of $d^{x}$. 
\label{imprk}
\end{rk}

For  the following axiom to make sense we impose a technical condition on the co-domains $\displaystyle V_{\varepsilon}(x)$: for any compact set $K \subset X$ there are $R=R(K) > 0$ and 
$\displaystyle \varepsilon_{0}= \varepsilon(K) \in (0,1)$  such that  
for all $\displaystyle u,v \in \bar{B}_{d}(x,R)$ and all $\displaystyle \varepsilon \in \Gamma$, $\displaystyle  \nu(\varepsilon) \in (0,\varepsilon_{0})$,  we have 
$$\delta_{\varepsilon}^{x} v \in W_{\varepsilon^{-1}}( \delta^{x}_{\varepsilon}u) \ .$$

With this assumption the following notation makes sense:
$$\Delta^{x}_{\varepsilon}(u,v) = \delta_{\varepsilon^{-1}}^{\delta^{x}_{\varepsilon} u} \delta^{x}_{\varepsilon} v . $$
The next axiom can now be stated: 
\begin{enumerate}
\item[{\bf A4.}] We have the limit 
$$\lim_{\varepsilon \rightarrow 0}  \Delta^{x}_{\varepsilon}(u,v) =  \Delta^{x}(u, v)  $$
uniformly with respect to $x, u, v$ in compact set. 
\end{enumerate}

\begin{defi}
A triple $(X,d,\delta)$ which satisfies A0, A1, A2, A3, but $\displaystyle d^{x}$ is degenerate for some 
$x\in X$, is called degenerate dilatation structure. 

If the triple $(X,d,\delta)$ satisfies A0, A1, A2, A3 and 
 $\displaystyle d^{x}$ is non-degenerate for any $x\in X$, then we call it  a weak dilatation structure. 

 If a weak dilatation structure satisfies A4 then we call it dilatation structure. 
 \end{defi}

\subsection{Groups with dilatations. Conical groups}

Metric tangent spaces sometimes have a group structure which is compatible 
with dilatations. This structure, of a group with dilatations, is interesting 
by itself. The notion has been introduced in \cite{buliga2}; we describe  it 
further.

Let $G$ be a topological group endowed with an uniformity such that the 
operation is uniformly continuous.  The  following description 
  is slightly non canonical, but is nevertheless motivated by the case of a 
Lie group endowed with a Carnot-Caratheodory  
distance induced by a left invariant distribution (see for example \cite{srlie1}, 
\cite{buliga2}).

We introduce first the double of $G$, as the group $G^{(2)} \ = \ G \times G$ 
with operation
$$(x,u) (y,v) \ = \ (xy, y^{-1}uyv)$$
The operation on the group $G$, seen as the function
$$op: G^{(2)} \rightarrow G \ , \ \ op(x,y) \ = \ xy$$
is a group morphism. Also the inclusions:
$$i': G \rightarrow G^{(2)} \ , \ \ i'(x) \ = \ (x,e) $$
$$i": G \rightarrow G^{(2)} \ , \ \ i"(x) \ = \ (x,x^{-1}) $$
are group morphisms.

\begin{defi}
\begin{enumerate}
\item[1.]
$G$ is an uniform group if we have two uniformity structures, on $G$ and
$G\times G$,  such that $op$, $i'$, $i"$ are uniformly continuous.

\item[2.] A local action of a uniform group $G$ on a uniform  pointed space $(X, x_{0})$ is a function
$\phi \in W \in \mathcal{V}(e)  \mapsto \hat{\phi}: U_{\phi} \in \mathcal{V}(x_{0}) \rightarrow
V_{\phi}  \in \mathcal{V}(x_{0})$ such that:
\begin{enumerate}
\item[(a)] the map $(\phi, x) \mapsto \hat{\phi}(x)$ is uniformly continuous from $G \times X$ (with product uniformity)
to  $X$,
\item[(b)] for any $\phi, \psi \in G$ there is $D \in \mathcal{V}(x_{0})$
such that for any $x \in D$ $\hat{\phi \psi^{-1}}(x)$ and $\hat{\phi}(\hat{\psi}^{-1}(x))$ make sense and   $\hat{\phi \psi^{-1}}(x) = \hat{\phi}(\hat{\psi}^{-1}(x))$.
\end{enumerate}

\item[3.] Finally, a local group is an uniform space $G$ with an operation defined in a neighbourhood of $(e,e) \subset G \times G$ which satisfies the uniform group axioms locally.
\end{enumerate}
\label{dunifg}
\end{defi}
Remark that a local group acts locally at left (and also by conjugation) on itself.

An uniform group, according to the definition \eqref{dunifg}, is a group $G$ such that left translations are uniformly continuous functions and the left action of $G$ on itself is uniformly continuous too. 
In order to precisely formulate this we need two uniformities: one on $G$ and another on $G \times G$. 

These uniformities should be compatible, which is achieved by saying that $i'$, $i"$ are uniformly continuous. The uniformity of the group operation is achieved by saying that the $op$ morphism is uniformly continuous.

\begin{defi}
A group with dilatations $(G,\delta)$ is a local uniform group $G$  with  a local action of $\Gamma$ (denoted by $\delta$), on $G$ such that
\begin{enumerate}
\item[H0.] the limit  $\displaystyle \lim_{\varepsilon \rightarrow 0} \delta_{\varepsilon} x \ = \ e$ exists and is uniform with respect to $x$ in a compact neighbourhood of the identity $e$.
\item[H1.] the limit
$$\beta(x,y) \ = \ \lim_{\varepsilon \rightarrow 0} \delta_{\varepsilon}^{-1}
\left((\delta_{\varepsilon}x) (\delta_{\varepsilon}y ) \right)$$
is well defined in a compact neighbourhood of $e$ and the limit is uniform.
\item[H2.] the following relation holds
$$ \lim_{\varepsilon \rightarrow 0} \delta_{\varepsilon}^{-1}
\left( ( \delta_{\varepsilon}x)^{-1}\right) \ = \ x^{-1}$$
where the limit from the left hand side exists in a neighbourhood of $e$ and is uniform with respect to $x$.
\end{enumerate}
\label{defgwd}
\end{defi}

These axioms are in fact a particular version of the axioms for a dilatation structure. We shall explain this a bit later.  

 Further we define conical local uniform groups.

\begin{defi}
A conical group $N$ is a local group with a local action of
$\Gamma$ by morphisms $\delta_{\varepsilon}$ such that
$\displaystyle \lim_{\varepsilon \rightarrow 0} \delta_{\varepsilon} x \ = \ e$ for any
$x$ in a neighbourhood of the neutral element $e$.
\end{defi}

The next proposition  explains why a conical group is the infinitesimal version of a group with 
dilatations.

\begin{prop}
Under the hypotheses H0, H1, H2 $(G,\beta, \delta)$ is a conical group, with operation 
$\beta$ and dilatations $\delta$.
\label{here3.4}
\end{prop}

Any group with dilatations has an associated dilatation structure on it.  In a group with dilatations $(G, \delta)$  we define dilatations based in any point $x \in G$ by 
 \begin{equation}
 \delta^{x}_{\varepsilon} u = x \delta_{\varepsilon} ( x^{-1}u)  . 
 \label{dilat}
 \end{equation}
 
\begin{defi} A normed group with dilatations $(G, \delta, \| \cdot \|)$ is a 
group with dilatations  $(G, \delta)$ endowed with a continuous norm  
function $\displaystyle \|\cdot \| : G \rightarrow \mathbb{R}$ which satisfies 
(locally, in a neighbourhood  of the neutral element $e$) the properties: 
 \begin{enumerate}
 \item[(a)] for any $x$ we have $\| x\| \geq 0$; if $\| x\| = 0$ then $x=e$, 
 \item[(b)] for any $x,y$ we have $\|xy\| \leq \|x\| + \|y\|$, 
 \item[(c)] for any $x$ we have $\displaystyle \| x^{-1}\| = \|x\|$, 
 \item[(d)] the limit 
$\displaystyle \lim_{\varepsilon \rightarrow 0} \frac{1}{\nu(\varepsilon)} \| \delta_{\varepsilon} x \| = \| x\|^{N}$ 
 exists, is uniform with respect to $x$ in compact set, 
 \item[(e)] if $\displaystyle \| x\|^{N} = 0$ then $x=e$.
  \end{enumerate}
  \label{dnco}
  \end{defi}
  
  It is easy to see that if $(G, \delta, \| \cdot \|)$ is a normed group with dilatations then $\displaystyle (G, \beta, \delta, \|\cdot\|^{N})$ is a normed conical group. The norm $\displaystyle \|\cdot\|^{N}$ satisfies the 
  stronger form of property (d) definition \ref{dnco}: for any $\varepsilon >0$ 
  $$ \| \delta_{\varepsilon} x \|^{N} = \varepsilon \| x \|^{N} .$$
  
 Normed conical groups generalize the notion of Carnot groups. 

In a normed group with dilatations we have a natural left invariant distance given by
\begin{equation}
d(x,y) = \| x^{-1}y\| . 
\label{dnormed}
\end{equation}

\begin{thm}
Let $(G, \delta, \| \cdot \|)$ be  a locally compact  normed group with dilatations. Then $(G, \delta, d)$ is 
a dilatation structure, where $\delta$ are the dilatations defined by (\ref{dilat}) and the distance $d$ is induced by the norm as in (\ref{dnormed}). 
\label{tgrd}
\end{thm}

\subsection{Tangent bundle of a dilatation structure}

\begin{thm}
Let $(X,d,\delta)$ be a weak dilatation structure. Then
\begin{enumerate}
\item[(a)] for all $x\in X$, $u,v \in X$ such that $\displaystyle d(x,u)\leq 1$ and $\displaystyle d(x,v) \leq 1$  and all $\mu \in (0,A)$ we have: 
$$d^{x}(u,v) \ = \ \frac{1}{\mu} d^{x}(\delta_{\mu}^{x} u , \delta^{x}_{\mu} v) \ .$$
We shall say that $d^{x}$ has the cone property with respect to dilatations. 
\item[(b)] we have  the following limit: 
$$\lim_{\varepsilon \rightarrow 0} \ \frac{1}{\varepsilon} \sup \left\{  \mid d(u,v) - d^{x}(u,v) \mid \mbox{ : } d(x,u) \leq \varepsilon \ , \ d(x,v) \leq \varepsilon \right\} \ = \ 0 \ .$$
Therefore  $(X,d)$ admits a metric tangent space 
at $x$, for any point $x\in X$. 
\end{enumerate}
\label{thcone}
\end{thm}

For the next theorem we need the previously introduced notion  of a conical (local) group.

\begin{thm}
Let $(X,d,\delta)$ be a dilatation structure. Then for any $x \in X$ the triple  $\displaystyle (U(x), \Sigma^{x}, \delta^{x})$ is a conical group. Moreover, left translations of this group are $\displaystyle d^{x}$ isometries. 
\label{tgene}
\end{thm}

The conical group $\displaystyle (U(x), \Sigma^{x}, \delta^{x})$ can be regarded as the tangent space 
of $(X,d, \delta)$ at $x$. Further will be denoted by: 
$\displaystyle T_{x} X =  (U(x), \Sigma^{x}, \delta^{x})$.

The following definition will be used in several further places. 

\begin{defi}
Let $(X,\delta,d)$ be a dilatation structure and $x\in X$ a point. 
In a neighbourhood $U(x)$ of $x$, for  any $\mu\in (0,1)$ 
we defined the distances:
$$(\delta^{x},\mu)(u,v) = \frac{1}{\mu} d(\delta^{x}_{\mu} u , \delta^{x}_{\mu} v) . $$
\label{drelative}
\end{defi}

\subsection{Topological considerations}
 
 In this subsection we compare various topologies and uniformities related to a dilatation structure.  
  
 The axiom A3 implies that for any $x \in X$ the function $\displaystyle d^{x}$ is continuous, therefore 
 open sets with respect to $\displaystyle d^{x}$ are open with respect to $d$. 
 
 If $(X,d)$ is separable and $\displaystyle d^{x}$ is non degenerate then $\displaystyle (U(x), d^{x})$ is also separable and the topologies of $d$ and $\displaystyle d^{x}$ are the same. Therefore $\displaystyle (U(x), d^{x})$ is also locally compact (and a set is compact  with respect to $d^{x}$ if and only if it is compact with respect to $d$). 
 
 If  $(X,d)$ is separable and $\displaystyle d^{x}$ is non degenerate then the uniformities induced by 
 $d$ and  $\displaystyle d^{x}$ are the same. Indeed, let 
 $$\left\{u_{n} \mbox{ : } n \in \mathbb{N}\right\} $$
 be a dense set in $U(x)$, with $\displaystyle x_{0}=x$. 
 We can embed $\displaystyle (U(x), 
 (\delta^{x}, \varepsilon))$ isometrically in a separable Banach space, for any $\varepsilon \in (0,1)$, by the function 
 $$\phi_{\varepsilon}(u) = \left( \frac{1}{\varepsilon} d(\delta^{x}_{\varepsilon}u,  \delta^{x}_{\varepsilon}x_{n}) - \frac{1}{\varepsilon} d(\delta^{x}_{\varepsilon}x,  \delta^{x}_{\varepsilon}x_{n})\right)_{n}  . $$
 A reformulation of point (a) in theorem \ref{thcone} is that on compact sets $\displaystyle \phi_{\varepsilon}$ uniformly converges to the isometric embedding of $\displaystyle (U(x), d^{x})$ 
 $$\phi(u) = \left(  d^{x}(u,  x_{n}) - d^{x}(x, x_{n})\right)_{n}  . $$
Remark that the uniformity induced by $(\delta,\varepsilon)$ is the same as the uniformity 
induced by $d$, and that it is the same induced from the uniformity on the separable Banach space by 
the embedding $\displaystyle \phi_{\varepsilon}$. We proved that the uniformities induced by 
 $d$ and  $\displaystyle d^{x}$ are the same.

\subsection{Equivalent dilatation structures}

\begin{defi}
Two dilatation structures $(X, \delta , d)$ and $(X, \overline{\delta} , \overline{d})$  are equivalent  if 
\begin{enumerate}
\item[(a)] the identity  map $\displaystyle id: (X, d) \rightarrow (X, \overline{d})$ is bilipschitz and 
\item[(b)]  for any $x \in X$ there are functions $\displaystyle P^{x}, Q^{x}$ (defined for $u \in X$ sufficiently close to $x$) such that  
\begin{equation}
\lim_{\varepsilon \rightarrow 0} \frac{1}{\varepsilon} \overline{d} \left( \delta^{x}_{\varepsilon} u ,  \overline{\delta}^{x}_{\varepsilon} Q^{x} (u) \right)  = 0 , 
\label{dequiva}
\end{equation}
\begin{equation}
 \lim_{\varepsilon \rightarrow 0} \frac{1}{\varepsilon} d \left( \overline{\delta}^{x}_{\varepsilon} u ,  
 \delta^{x}_{\varepsilon} P^{x} (u) \right)  = 0 , 
\label{dequivb}
\end{equation}
uniformly with respect to $x$, $u$ in compact sets. 
\end{enumerate}
\label{dilequi}
\end{defi}

\begin{prop}
Two dilatation structures $(X, \delta , d)$ and $(X, \overline{\delta} , \overline{d})$  are equivalent  if and 
only if 
\begin{enumerate}
\item[(a)] the identity  map $\displaystyle id: (X, d) \rightarrow (X, \overline{d})$ is bilipschitz and 
\item[(b)]  for any $x \in X$ there are functions $\displaystyle P^{x}, Q^{x}$ (defined for $u \in X$ sufficiently close to $x$) such that  
\begin{equation}
\lim_{\varepsilon \rightarrow 0}  \left(\overline{\delta}^{x}_{\varepsilon}\right)^{-1}  \delta^{x}_{\varepsilon} (u) = Q^{x}(u) , 
\label{dequivap}
\end{equation}
\begin{equation}
 \lim_{\varepsilon \rightarrow 0}  \left(\delta^{x}_{\varepsilon}\right)^{-1}  \overline{\delta}^{x}_{\varepsilon} (u) = P^{x}(u) , 
\label{dequivbp}
\end{equation}
uniformly with respect to $x$, $u$ in compact sets. 
\end{enumerate}
\label{pdilequi}
\end{prop}

 The next theorem shows a link between the tangent bundles of equivalent dilatation structures. 
 
 \begin{thm} 
 Let $(X, \delta , d)$ and $(X, \overline{\delta} , \overline{d})$  be  equivalent dilatation structures. Suppose that for any $x \in X$ the distance $d^{x}$ is non degenerate. Then for any $x \in X$ and 
 any $u,v \in X$ sufficiently close to $x$ we have:
 \begin{equation}
 \overline{\Sigma}^{x}(u,v) = Q^{x} \left( \Sigma^{x} \left( P^{x}(u) , P^{x}(v) \right)\right) . 
 \label{isoequiv}
 \end{equation}
 The two tangent bundles  are therefore isomorphic in a natural sense. 
 \label{tisoequiv}
 \end{thm}

\subsection{Differentiable functions}
 
 Dilatation structures allow to define differentiable functions. The idea is to 
keep only  one relation from definition \ref{dilequi}, namely (\ref{dequiva}). 
We also renounce to uniform convergence with respect 
 to $x$ and $u$, and we replace this with uniform convergence in the "$u$" variable,  
with a conical group morphism condition for the derivative.

\begin{defi}
 Let $(N,\delta)$ and $(M,\bar{\delta})$ be two conical groups. A function $f:N\rightarrow M$ is a conical group morphism if $f$ is a group morphism and for any $\varepsilon>0$ and $u\in N$ we have 
 $\displaystyle f(\delta_{\varepsilon} u) = \bar{\delta}_{\varepsilon} f(u)$. 
\label{defmorph}
\end{defi}

 \begin{defi}
 Let $(X, \delta , d)$ and $(Y, \overline{\delta} , \overline{d})$ be two dilatation structures and $f:X \rightarrow Y$ be a continuous function. The function $f$ is differentiable in $x$ if there exists a 
 conical group morphism  $\displaystyle Q^{x}:T_{x}X\rightarrow T_{f(x)}Y$, defined on a neighbourhood of $x$ with values in  a neighbourhood  of $f(x)$ such that 
\begin{equation}
\lim_{\varepsilon \rightarrow 0} \sup \left\{  \frac{1}{\varepsilon} \overline{d} \left( f\left( \delta^{x}_{\varepsilon} u\right) ,  \overline{\delta}^{f(x)}_{\varepsilon} Q^{x} (u) \right) \mbox{ : } d(x,u) \leq \varepsilon \right\}   = 0 , 
\label{edefdif}
\end{equation}
The morphism $\displaystyle Q^{x}$ is called the derivative of $f$ at $x$ and will be sometimes denoted by $Df(x)$.

The function $f$ is uniformly differentiable if it is differentiable everywhere and the limit in (\ref{edefdif}) 
is uniform in $x$ in compact sets. 
\label{defdif}
\end{defi}
 
 A trivial way to obtain a differentiable function (everywhere)  is to modify the dilatation structure on the target space. 
 
 \begin{defi}
 Let  $(X, \delta , d)$ be a  dilatation structure and $f:(X, d) \rightarrow (Y, \overline{d})$ be a bilipschitz  and surjective  function. We define then the transport of $(X, \delta , d)$ by $f$, named $(Y, f*\delta , \overline{d})$, by: 
 $$\left( f*\delta\right)^{f(x)}_{\varepsilon} f(u) = f \left( \delta^{x}_{\varepsilon} u \right) . $$
 \label{ddif}
 \end{defi}

 The relation of differentiability with equivalent dilatation structures is given by the following simple proposition.

 \begin{prop}
 Let  $(X, \delta , d)$ and $(X, \overline{\delta} , \overline{d})$ be two dilatation structures and $f:(X, d) \rightarrow (X, \overline{d})$ be a bilipschitz  and surjective  function. The dilatation structures $(X, \overline{\delta} , \overline{d})$ and $(X, f*\delta , \overline{d})$ are equivalent if and only if $f$ and $\displaystyle f^{-1}$ are uniformly  differentiable. 
 \label{peqd}
 \end{prop}

 We shall prove now the chain rule for derivatives, after we elaborate a bit over the definition \ref{defdif}.

 Let $(X, \delta , d)$ and $(Y, \overline{\delta} , \overline{d})$ be two dilatation structures and $f:X \rightarrow Y$  a function differentiable in $x$. The derivative of $f$ in $x$ is a 
 conical group morphism  $\displaystyle Df(x):T_{x}X\rightarrow T_{f(x)}Y$, which means that $Df(x)$ is 
 defined on a open set around $x$ with values in a open set around $f(x)$, having the properties:
 \begin{enumerate}
 \item[(a)] for any $u,v$ sufficiently close to $x$ 
 $$Df(x)\left(\Sigma^{x}(u,v)\right) = \Sigma^{f(x)}\left(Df(x)(u), Df(x)(v)\right)  , $$
 \item[(b)] for any $u$ sufficiently close to $x$ and any $\varepsilon \in (0,1]$ 
 $$Df(x)\left(\delta^{x}_{\varepsilon} u\right) = \bar{\delta}^{f(x)}_{\varepsilon}\left(Df(x)(u)\right) , $$
 \item[(c)] the function $Df(x)$ is continuous, as uniform limit of continuous 
functions. Indeed, the relation 
 (\ref{edefdif}) is equivalent to the existence of the uniform limit 
(with respect to $u$ in compact sets)
 $$Df(x)(u) = \lim_{\varepsilon\rightarrow 0} \bar{\delta}^{f(x)}_{\varepsilon^{-1}} \left( f\left( \delta_{\varepsilon}^{x} u\right)\right) . $$ 
 \end{enumerate}
 
 From (\ref{edefdif}) alone and axioms of dilatation structures we can prove properties (b) and (c). 
 We can reformulate therefore the definition of the derivative by asking that $Df(x)$ exists as an uniform 
 limit (as in point (c) above) and that $Df(x)$ has the property (a) above. 
 
 From these considerations the chain rule for derivatives is straightforward. 
 
 \begin{prop}
 Let $(X, \delta , d)$,  $(Y, \overline{\delta} , \overline{d})$ and $(Z, \hat{\delta} , \hat{d})$ be three  dilatation structures and $f:X \rightarrow Y$  a continuous function differentiable in $x$, $g:Y \rightarrow Z$  a continuous function differentiable in $f(x)$. Then $gf:X \rightarrow Z$ is   differentiable in $x$ and 
$$D gf (x) = Dg(f(x)) Df(x) . $$
\label{pchain}
\end{prop}

\paragraph{Proof.}
Use property (b) for proving that $Dg(f(x)) Df(x)$ satisfies (\ref{edefdif}) for the function $gf$ and $x$. 
Both $Dg(f(x))$ and $Df(x)$ are conical group morphisms,  therefore  $Dg(f(x)) Df(x)$ is a conical group 
morphism too. We deduce that  $Dg(f(x)) Df(x)$ is the derivative of $gf$ in $x$. \quad $\|$

\subsection{Some induced dilatation structures}
\label{induced}

\begin{prop}
For any $u,v\in U(x)$ let us  define 
$$\hat{\delta}_{\varepsilon}^{u} v = \Sigma^{x}_{\mu}(u, \delta^{\delta^{x}_{\mu}u}_{\varepsilon} 
\Delta^{x}_{\mu}(u,v)) = \delta^{x}_{\mu^{-1}} \delta_{\varepsilon}^{\delta_{\mu}^{x} u} \delta_{\mu}^{x} v  . $$
Then $\displaystyle (U(x),\hat{\delta}, (\delta^{x},\mu))$ is a dilatation structure. 
\label{pshift}
\end{prop}

\paragraph{Proof.}
We have to check the axioms.  The first part of axiom A0 is an easy consequence of theorem 
\ref{thcone} for $(X,\delta,d)$. The second part of A0, A1 and A2 are true based on simple computations. 

The first interesting fact is related to axiom A3. Let us compute, for $v,w \in U(x)$, 
$$\frac{1}{\varepsilon} (\delta^{x},\mu)(\hat{\delta}^{u}_{\varepsilon} v, \hat{\delta}^{u}_{\varepsilon} w) = 
\frac{1}{\varepsilon \mu} d( \delta^{x}_{\mu} \hat{\delta}^{u}_{\varepsilon} v, \delta^{x}_{\mu}
\hat{\delta}^{u}_{\varepsilon} w) = $$
$$ = \frac{1}{\varepsilon \mu} d( \delta_{\varepsilon}^{\delta^{x}_{\mu} u} \delta_{\mu}^{x} v ,  
\delta_{\varepsilon}^{\delta^{x}_{\mu} u} \delta_{\mu}^{x} w) = \frac{1}{\varepsilon\mu} d(  
\delta_{\varepsilon \mu}^{\delta^{x}_{\mu} u} \Delta_{\mu}^{x}(u,v),  \delta_{\varepsilon \mu}^{\delta^{x}_{\mu} u} \Delta_{\mu}^{x}(u,w)) =  $$
$$= (\delta^{\delta_{\mu}^{x} u}, \varepsilon \mu) (  \Delta_{\mu}^{x}(u,v) ,  \Delta_{\mu}^{x}(u,w)) . $$
The axiom A3 is then a consequence of axiom A3 for $(X,\delta,d)$ and we have 
$$\lim_{\varepsilon\rightarrow 0} \frac{1}{\varepsilon} (\delta^{x},\mu)(\hat{\delta}^{u}_{\varepsilon} v, \hat{\delta}^{u}_{\varepsilon} w) = d^{\delta_{\mu}^{x} u} (  \Delta_{\mu}^{x}(u,v) ,  \Delta_{\mu}^{x}(u,w)) . $$
The axiom A4 is also a straightforward consequence of A4 
for $(X,\delta,d)$ and is left to the reader. \quad $\square$

 The proof of the following proposition is an easy computation, of the same 
type as in the lines above, therefore we shall not write it here. 

\begin{prop}
With the same notations as in proposition \ref{pshift}, the transformation $\displaystyle \Sigma^{x}_{\mu}(u, \cdot)$ is an isometry from $\displaystyle (\delta^{\delta^{x}_{\mu} u}, \mu)$ to $\displaystyle (\delta^{x}, \mu)$. Moreover, we have 
$$ \Sigma^{x}_{\mu} (u, \delta_{\mu}^{x} u) = u . $$
\label{isoshift}
\end{prop}

These two propositions show  that on a dilatation structure we almost have translations (the infinitesimal sums), which are almost isometries (that is, not with respect to the distance $d$, but with respect to distances of type $\displaystyle (\delta^{x},\mu)$). It is almost as if we were working with  a 
conical group, only that we have to use families of distances and to make small shifts in the tangent space (as in the last formula in the proof of propositon \ref{pshift}). Moreover, in a very precise way everything converges as $\mu\rightarrow 0$ to the right thing. 

\section{The linear group of a dilatation structure}

\begin{defi}
Let $(X,d,\delta)$ be a weak dilatation structure. A transformation $A:X\rightarrow X$ is linear if it 
is Lipschitz and it commutes with dilatations in the following sense: for any $x\in X$, $u \in U(x)$ and 
$\varepsilon \in \Gamma$, $\nu(\varepsilon) < 1$, if  $A(u) \in U(A(x))$ then  
$$ A \delta^{x}_{\varepsilon} = \delta^{A(x)} A(u) \quad  .$$
The  group of linear transformations, denoted by $GL(X,d,\delta)$ is formed by all invertible and 
bi-lipschitz linear transformations of $X$. 
\label{defgl}
\end{defi}

$GL(X,d,\delta)$ is  a (local) group. Indeed, we start from the  
remark that if $A$ is Lipschitz then there exists $C>0$ such that for all $x\in X$ and $u \in B(x,C)$ we have $A(u)\in U(A(x))$.  The inverse of $A \in GL(X,d,\delta)$ is then linear. Same considerations apply for the composition of two linear, bi-lipschitz and invertible transformations. 

In the particular case of  example \ref{401}, namely  $X$ finite dimensional real, normed vector space, 
$d$ the distance given by the norm, $\Gamma = (0,+\infty)$ and dilatations 
$$\delta_{\varepsilon}^{x} u = x + \varepsilon(u-x) \quad , $$
a linear transformations in the sense of definition \ref{defgl} is an affine transformation of the vector 
space $X$. 

Linear transformations have nice properties which justify the name "linear".  We shall use further the  
(sum and difference) operations 
$$\Sigma_{\varepsilon}^{x}(u,v) = \left( \delta_{\varepsilon}^{x}\right)^{-1} \, \delta^{\delta^{x}_{\varepsilon} u}_{\varepsilon} v \quad , \quad  \Delta_{\varepsilon}^{x}(u,v) = \left( \delta_{\varepsilon}^{\delta^{x}_{\varepsilon}}\right)^{-1} \, \delta^{x}_{\varepsilon} v $$
and the inverse function $\displaystyle inv^{x}_{\varepsilon} u \, = \Delta^{x}_{\varepsilon}(u,x)$.

\begin{prop}
Let $(X,d,\delta)$ be a weak dilatation structure and $A:X\rightarrow X$ a linear transformation. Then: 
\begin{enumerate}
\item[(a)] for all $x\in X$, $u,v \in U(x)$ sufficiently close to $x$, we have: 
$$A \, \Sigma_{\varepsilon}^{x}(u,v)  = \Sigma_{\varepsilon}^{A(x)}(A(u),A(v)) \quad . $$
\item[(b)] or all $x\in X$, $u \in U(x)$ sufficiently close to $x$, we have: 
$$ A \, inv^{x}(u) = \, inv^{A(x)} A(u) \quad . $$
\item[(c)]  for all $x\in X$ the transformation $A$ is derivable and the derivative equals $A$. 
\end{enumerate}
\label{plinear}
\end{prop}

\paragraph{Proof.}
Straightforward, just use the commutation with dilatations. $\quad \square$

This is important because  the sum, difference, inverse operations induced by a dilatation structure 
give to the space $X$ almost the structure of an affine space. We collect some results from 
\cite{buligadil1} section 4.2 , regarding the properties of these operations. Only the last point is new, 
but with straightforward proof. 

\begin{thm}
Let $(X,d,\delta)$ be a weak dilatation structure. Then, for any $x \in X$, $\varepsilon \in \Gamma$, 
$\nu(\varepsilon) < 1$,     we have: 
\begin{enumerate}
\item[(a)] for any $u\in U(x)$,  $\displaystyle \Sigma_{\varepsilon}^{x}(x,u) = u $  .
\item[(b)] for any $u\in U(x)$ the functions $\displaystyle  \Sigma_{\varepsilon}^{x}(u,\cdot)$ and 
$\displaystyle  \Delta_{\varepsilon}^{x}(u, \cdot)$ are inverse one to another. 
\item[(c)] the inverse function is shifted involutive:  for any $u\in U(x)$, 
$$inv^{\delta_{\varepsilon}^{x} u}_{\varepsilon} \, inv^{x}_{\varepsilon} (u) = u \quad . $$
\item[(d)] the sum operation is shifted associative: for  any $u, v, w$ sufficiently close to $x$ we have 
$$\Sigma_{\varepsilon}^{x} \left( u, \Sigma_{\varepsilon}^{\delta^{x}_{\varepsilon} u} (v, w)\right) = 
\Sigma_{\varepsilon}^{x} ( \Sigma^{x}(u,v), w) \quad . $$
\item[(e)] the difference, inverse and sum operations are related by 
$$ \Delta_{\varepsilon}^{x}(u,v) = \Sigma_{\varepsilon}^{\delta_{\varepsilon}^{x} u}
 \left( inv_{\varepsilon}^{x}(u), v\right) \quad , $$
 for any $u,v$ sufficiently close to $x$. 
 \item[(f)] for any $u,v$ sufficiently close to $x$ and $\mu \in \Gamma$, 
$\nu(\mu) < 1$,     we have: 
$$\Delta^{x}_{\varepsilon} \left( \delta^{x}_{\mu} u,  \delta^{x}_{\mu} v \right) =  \delta^{ \delta^{x}_{\epsilon \mu} u}_{\mu} \Delta^{x}_{\varepsilon \mu} (u,v) \quad . $$
\end{enumerate}
\label{opcollection}
\end{thm}

Remark that in principle the "translations" $\displaystyle \Sigma^{x}_{\varepsilon}(u, \cdot)$  are not linear. Nevertheless, they commute with dilatation in a known way, according to point (f) theorem  
\ref{opcollection}. This is important, because the transformations $\displaystyle \Sigma^{x}_{\varepsilon}(u, \cdot)$ really behave as translations, as explained in subsection \ref{induced}. 

The reason for which translations are not linear is that dilatations are not linear. In the case of strong 
dilatation structures, this happens only when we are in a conical group. 

\begin{thm}
Let $(X,d,\delta)$ be a weak dilatation structure. 
\begin{enumerate}
\item[(a)] If dilatations are linear then all    transformations 
$\displaystyle \Delta^{x}_{\varepsilon}(u, \cdot)$ are linear for any $u \in X$. 
\item[(b)]   If the dilatation structure is strong then dilatations are linear if and only if the dilatations  come  from the  dilatation structure of a conical group. 
\end{enumerate}
\label{tdilatlin}
\end{thm}

\paragraph{Proof.}
(a) If dilatations are linear, then let $\varepsilon, \mu \in \Gamma$, $\nu(\varepsilon), 
\nu(\mu) \leq 1$, and $x, y, u, v \in X$  such that the following computations make sense. We have: 
$$\Delta^{x}_{\varepsilon} (u, \delta_{\mu}^{y} v ) = \delta_{\varepsilon^{-1}}^{\delta^{x}_{\varepsilon} u} \delta_{\varepsilon}^{x} \delta_{\mu}^{y} v \quad . $$
Let $\displaystyle A_{\varepsilon} = \delta_{\varepsilon^{-1}}^{\delta^{x}_{\varepsilon} u}$. We compute: 
$$\delta_{\mu}^{\Delta_{\varepsilon}^{x}(u,y)} \Delta_{\varepsilon}^{x} (u,v)  = \delta_{\mu}^{A_{\varepsilon} \delta_{\varepsilon}^{x} y} A_{\varepsilon} \delta_{\varepsilon}^{x} v \quad . $$
We use twice the linearity of dilatations:  
$$ \delta_{\mu}^{\Delta_{\varepsilon}^{x}(u,y)} \Delta_{\varepsilon}^{x} (u,v)  = A_{\varepsilon} 
\delta_{\mu}^{\delta_{\varepsilon}^{x} y} \delta_{\varepsilon}^{x} v = 
 \delta_{\varepsilon^{-1}}^{\delta^{x}_{\varepsilon} u} \delta_{\varepsilon}^{x} \delta_{\mu}^{y} v \quad . $$
We proved that: 
$$\Delta^{x}_{\varepsilon} (u, \delta_{\mu}^{y} v ) =  \delta_{\mu}^{\Delta_{\varepsilon}^{x}(u,y)} \Delta_{\varepsilon}^{x} (u,v)  \quad , $$
which is the conclusion of   the part (a).

(b) Suppose that the dilatation structure is strong. If dilatations are linear, then by point (a) the  transformations $\displaystyle \Delta^{x}_{\varepsilon}(u, \cdot)$$\delta$ are linear for any $u \in X$. Then, with notations made before, for $y = u$ we get 
$$\Delta^{x}_{\varepsilon} (u, \delta_{\mu}^{u} v ) =  \delta_{\mu}^{\delta_{\varepsilon}^{x} u} \Delta_{\varepsilon}^{x} (u,v)  \quad , $$
which implies 
$$\delta_{\mu}^{u} v = \Sigma^{x}_{\varepsilon} ( u, \delta^{x}_{\mu} \Delta_{\varepsilon}^{x}(u,v)) \quad . $$
We pass to the limit with $\varepsilon \rightarrow 0$ and we obtain: 
$$ \delta_{\mu}^{u} v = \Sigma^{x}( u, \delta^{x}_{\mu} \Delta^{x}(u,v)) \quad . $$
We recognize at the  right hand side the dilatations associated to the conical group 
$\displaystyle T_{x} X$. 

The opposite implication is straightforward, because the dilatation structure of any conical group 
is linear. 
\quad $\square$

  \section{First examples of dilatation structures}

In this section we give several examples of dilatation structures, 
which share some common features.

\begin{exemplu}
The first example is known to everybody: take $\displaystyle (X,d)  =  ( \mathbb{R}^{n}, d_{E})$, with usual (euclidean) dilatations $\displaystyle \delta^{x}_{\varepsilon}$, that  is: 
$$d_{E}(x,y) = \| x-y \| \ \ \ , \ \ \delta_{\varepsilon}^{x} y \ = \ x + \varepsilon (y-x) \ .$$
Dilatations are defined everywhere. The group $\Gamma$ is $(0,+\infty)$ and the function $\nu$ is the identity. 

There are few things to check:  axioms 0,1,2 are obviously true. For axiom A3, remark that 
for any $\varepsilon > 0$, $x,u,v \in X$ we have: 
$$\frac{1}{\varepsilon} d_{E}(\delta^{x}_{\varepsilon} u , \delta^{x}_{\varepsilon} v ) \ = \ d_{E}(u,v) \ , $$
therefore for any $x \in X$ we have $\displaystyle d^{x} = d_{E}$. 

Finally, let us check the axiom A4. For any $\varepsilon > 0$ and $x,u,v \in X$ we have
$$\delta_{\varepsilon^{-1}}^{\delta_{\varepsilon}^{x} u} \delta_{\varepsilon}^{x} v \ = \ 
x + \varepsilon  (u-x) + \frac{1}{\varepsilon} \left( x+ \varepsilon(v-x) - x - \varepsilon(u-x) \right) \ = \ $$
$$ = \ x + \varepsilon  (u-x) + v - u$$ 
therefore this quantity converges to 
$$x + v - u \ = \ x + (v - x) - (u - x)$$
as $\varepsilon \rightarrow 0$. The axiom A4 is verified. 
\label{401}
\end{exemplu}

\subsection{Standard dilatation structures}

\begin{exemplu} Take now $\phi: \mathbb{R}^{n} \rightarrow \mathbb{R}^{n}$ a bi-Lipschitz diffeomorphism. Then 
we can define the dilatation structure: $X \ = \ \mathbb{R}^{n}$, 
$$d_{\phi}(x,y) \ = \ \| \phi(x) - \phi(y) \| \ \ \ , \ \ \delta_{\varepsilon}^{x} y \ = \ x+ \varepsilon(y-x) \ ,$$
or the equivalent dilatation structure: 
$X \ = \ \mathbb{R}^{n}$, 
$$d_{\phi}(x,y) \ = \ \| x - y \| \ \ \ , \ \ \delta_{\varepsilon}^{x} y \ = \ \phi^{-1} \left( \phi(x) + \varepsilon (\phi(y) - \phi(x))\right) .$$
In this example (look at its first  version) the distance $\displaystyle d_{\phi}$ is not equal to $\displaystyle d^{x}$. Indeed, a direct calculation shows that 
$$d^{x}(u,v) \ = \ \| D\phi(x) (v-u) \| \ .$$
The axiom A4 gives the same result as previously. 
\label{411}
\end{exemplu}

\begin{exemplu}
Because dilatation structures are defined by local requirements, we can easily define dilatation 
structures on riemannian manifolds, using particular atlases of the manifold and the riemannian 
distance (infimum of length of curves joining two points). Note that any finite dimensional manifold can be endowed with a riemannian metric. This class of examples covers all dilatation structures used in 
differential geometry. The axiom A4 gives an operation of addition of vectors 
in the tangent space  (compare with Bella\"{\i}che \cite{bell} last section).  
\label{411riemann}
\end{exemplu}

There is a version of the snowflake construction for dilatation structures. 
It is stated in the next proposition, ehich has a  straightforward proof. 

\begin{prop}
If $(X,d,\delta)$ is a dilatation structure then $(X,d_{a}, \delta(a))$ is
 also a dilatation structure, for any $a \in (0,1]$, where 
$$d_{a}(x,y) \ = \  d(x,y)^{a} \ \ , \ \ \delta(a)_{\varepsilon}^{x}\ = \ \delta^{x}_{\varepsilon^{\frac{1}{a}}} \ .$$
\label{psnow}
\end{prop}

\begin{exemplu} In particular we get a  snowflake variation of the euclidean 
case: $X =  \mathbb{R}^{n}$ and  for any $a \in (0,1]$ take 
$$ d_{a}(x,y) \ = \ \| x-y \|^{\alpha} \ \ \ ,  \ \ \delta^{x}_{\varepsilon} y \ =  \ x + \varepsilon^{\frac{1}{a}} (y - x) \ .$$
\label{421}
\end{exemplu}

\subsection{Nonstandard dilatations in the euclidean space}

\begin{exemplu} Take $\displaystyle X = \mathbb{R}^{2}$ with the euclidean distance. For any $z \in \mathbb{C}$ of the 
form $z= 1+ i \theta$ we define dilatations 
$$\delta_{\varepsilon} x = \varepsilon^{z} x  \ .$$
It is easy to check that $(X,\delta, +, d)$ is a conical group, equivalenty that the dilatations 
$$\delta^{x}_{\varepsilon} y = x + \delta_{\varepsilon} (y-x)  \ .$$
form a linear dilatation structure with the euclidean distance. 

Two such dilatation structures (constructed with the help of complex numbers 
$1+ i \theta$ and $1+ i \theta'$) are equivalent if and only if $\theta = \theta'$.  

There are two other interesting  properties of these dilatation structures. 
The first is that if $\theta \not = 0$ then there are no non trivial Lipschitz curves in $X$ which are differentiable almost everywhere. 

The second property is that any holomorphic and Lipschitz function from $X$ to $X$ (holomorphic in the 
usual sense on $X = \mathbb{R}^{2} = \mathbb{C}$) is differentiable almost everywhere, but there are 
Lipschitz functions from $X$ to $X$ which are not differentiable almost everywhere (suffices to take a 
$\displaystyle \mathcal{C}^{\infty}$ function from  $\displaystyle \mathbb{R}^{2}$ to $\displaystyle \mathbb{R}^{2}$ which is not holomorphic). 
\label{exemp1}
\end{exemplu}

 Take now a one parameter group of linear transformations $\displaystyle s \in \mathbb{R} \mapsto  A_{s}$ in $\displaystyle X = \mathbb{R}^{2}$ such that $\displaystyle A_{s} \rightarrow 0$ as $s \rightarrow -\infty$. Such a group generates functions 
$$\delta_{\varepsilon} = A_{\log \varepsilon}  \ .$$
which can be used to construct dilatation structures. 

It is well known that, up to conjugation,  there are only three such 
one-parameter  groups in $\displaystyle \mathbb{R}^{2}$. 

\begin{exemplu}The first group generates diagonal functions $\delta$ with the form 
$$\delta_{\varepsilon} (x_{1}, x_{2}) = (\varepsilon^{\alpha} x_{1}, \varepsilon^{\beta} x_{2})  \ .$$
By the snowflake  construction we can find a distance on $X$ such that we get a dilatation structure. 
\label{exemp21}
\end{exemplu}

\begin{exemplu} The second group generates  functions $\delta$ with the form 
$$\delta_{\varepsilon} x = \varepsilon^{z} x  \ ,$$
with $Re  \                 z>0$. We can choose a distance on $X$ such that we have again a dilatation structure  (just combine the first example in this subsection with the snowflake construction). 
\label{exemp22}
\end{exemplu}

Finally, the third group generates functions $\delta$ with a different form. Modulo the snowflake 
construction the functions $\delta$ have the form: 
$$\delta_{\varepsilon} (x_{1}, x_{2}) = (\varepsilon x_{1} + \varepsilon \log (\varepsilon) x_{2}, \varepsilon x_{2})  \ .$$
If we choose the distance $d$ to be the euclidean distance then we verify all the axioms excepting the axiom A3.

Let now $\displaystyle \delta_{\varepsilon}$ be a one parameter group of linear invertible 
transformations (in multiplicative form) on $\displaystyle \mathbb{R}^{n}$, such that $\displaystyle \delta_{\varepsilon}$ converges to $0$ as $\varepsilon$ goes to $0$.  To any 
$\displaystyle x\in\mathbb{R}^{n}$ we associate, in a continuous way, 
 a linear invertible transformation of 
$\displaystyle \mathbb{R}^{n}$, denoted by $A(x)$. We define now 
$$\delta_{\varepsilon}^{x} y  = x + A(x) \delta_{\varepsilon}A(x)^{-1} ( y-x) . $$
We want to know if this is a dilatation structure on $(\mathbb{R}^{n},d)$, where $d$ is the euclidean distance. 

We have to check only axioms A3 and A4. For this notice that for any 
$\displaystyle u,v\in \mathbb{R}^{n}$ and $\varepsilon> 0$ 
$$\frac{1}{\varepsilon} d(\delta^{x}_{\varepsilon} u, \delta^{x}_{\varepsilon} v) = \| A(x) 
\frac{1}{\varepsilon} \delta_{\varepsilon} A(x)^{-1} (u-v)\| , $$
therefore the axiom A3 is satisfied if $\displaystyle \frac{1}{\varepsilon} \delta_{\varepsilon} $ 
converges to an invertible transformation. Assume that there is invertible (and linear, as limit 
of linear transformations) function $P$ such that 
$$\lim_{\varepsilon\rightarrow 0} \frac{1}{\varepsilon} \delta_{\varepsilon} = P . $$
In order to check the axiom A4 we compute: 
$$\Delta^{x}_{\varepsilon}(u,v) = x+A(\delta_{\varepsilon}^{x} u) 
\delta_{\varepsilon^{-1}}A(\delta_{\varepsilon}^{x} u)^{-1}A(x)\delta_{\varepsilon}A(x)^{-1}(v-x)  =  $$
$$= x+ \varepsilon A(\delta_{\varepsilon}^{x} u) 
\delta_{\varepsilon^{-1}}A(\delta_{\varepsilon}^{x} u)^{-1} \left( A(x)\frac{1}{\varepsilon}
\delta_{\varepsilon}A(x)^{-1}(v-x)\right) . $$ 
This expression converges if and only if the following limit exists: 
$$\lim_{\varepsilon\rightarrow 0} \varepsilon A(\delta_{\varepsilon}^{x} u) 
\delta_{\varepsilon^{-1}}A(\delta_{\varepsilon}^{x} u)^{-1}  . $$
The map  $x\mapsto A(x)$ continuous, threfore we  have the equivalence: $\displaystyle 
(\mathbb{R}^{n},d,\delta)$ is a dilatation structure if and only if we have the limit 
$$\lim_{\varepsilon\rightarrow 0} \varepsilon \delta_{\varepsilon^{-1}}A(\delta_{\varepsilon}^{x} u)^{-1}    . $$
This limit may not exist (or it may be infinite), thus providing examples of a structure which satisfies all axioms excepting A4. 

We  state as  a conclusion of this subsection:

\begin{thm}
The axioms A3 and A4 are independent of the rest of axioms. 
\label{tindep}
\end{thm}

\subsection{Ultrametric valued fields}

\begin{defi}
Let $R$ be a commutative ring with unity $1$. A function
 $\displaystyle N:R \rightarrow [0,+\infty)$ is a  norm  on 
$R$ if the following are true. 
\begin{enumerate}
\item[(a)] $N(x)=0$ if and only if $x=0$.
\item[(b)] For all $x,y\in R$ $N(xy) = N(x)N(y)$. 
\item[(c)] For all $x,y\in R$ $N(x+y)\leq N(x)+N(y)$.
\end{enumerate}
\label{defring}
\end{defi}
In case of a field, the set of   norms  of non zero elements of $R$ is called 
value group and denoted by  $\Theta \subset (0,+\infty)$.

Particular examples of  complete valued fields are  fields of $p$-adic numbers $\displaystyle 
\mathbb{Q}_{p}$ ($p$ prime). Good references are Schikhof \cite{schik}, or Bachman \cite{bachman}. 

 For such fields the value group $\Theta$ is discrete and condition 
(c) in definition above is strengthened to 
$$\forall x,y\in R \quad N(x+y)\leq \max\left\{ N(x), N(y)\right\} \quad .$$

\begin{exemplu}
Let $K$ be an ultrametric, complete, valued field, $$\displaystyle \bar{B}(0,1)= \left\{ x\in K \, \mid \, 
\mid x \mid \leq 1\right\}\quad , $$
 $$\displaystyle A= \left\{ x\in K \, \mid \, x\not = 0 , 
\mid x \mid \leq 1 , \mid 1-x\mid \leq 1 \right\}\quad . $$
Then $\displaystyle \bar{B}(0,1)$ is a compact set, a valued ring, $A$ is a semigroup with multiplication and 
$$\delta^{x}_{\varepsilon}: \bar{B}(0,1)\rightarrow \bar{B}(0,1) \quad , \quad \delta^{x}_{\varepsilon} y 
= x+ \varepsilon(y-x) \quad , $$
is well defined for any $\displaystyle x\in \bar{B}(0,1), \varepsilon\in A$. 

Moreover, the functions $\delta$ and the norm on $K$ define a dilatation structure 
$\displaystyle (\bar{B}(0,1), d, \delta)$.

Indeed,  the norm is ultrametric, therefore  $\displaystyle \bar{B}(0,1)$ is a compact set, a valued ring, and  $A$ is a semigroup with multiplication. We have to check that $\displaystyle \delta^{x}_{\varepsilon}$ is well defined  for any $\displaystyle x\in \bar{B}(0,1), \varepsilon\in A$. But this is straightforward: 
take any $\displaystyle y\in \bar{B}(0,1)$. Then 
$$\mid   \delta^{x}_{\varepsilon} y \mid \leq \max \left\{ \mid \varepsilon \mid \mid y\mid, \mid 1-\varepsilon
\mid \mid x\mid \right\} \leq 1 \quad . $$
The proof that we have here a dilatation structure is formally identical with example \ref{401}.
\label{tinteg}
\end{exemplu}

\section{Dilatation structures on the boundary of the dyadic tree}

Dilatation structures on the boundary of the dyadic tree will have a simpler form than general, mainly because the distance is ultrametric.

The boundary of the dyadic tree identifies with $\displaystyle X^{\omega}$, for 
$\displaystyle X = \left\{ 0,1 \right\}$, and also with the ring of dyadic integers $\displaystyle \mathbb{Z}_{2}$. Use shall use 
the usual, ultrametric distance, denoted by  $d$, on this set.

We  take the group $\Gamma$ to be the set of integer powers of $2$, seen as a 
subset of dyadic numbers. Thus for any $p \in \mathbb{Z}$ the element 
$\displaystyle 2^{p} \in \mathbb{Q}_{2}$ belongs to $\Gamma$. The operation is 
the multiplication of dyadic numbers and the morphism $\nu : \Gamma \rightarrow (0,+\infty)$ is defined by 
$$\nu(2^{p}) = d(0, 2^{p}) = \frac{1}{2^{p}} \in (0,+\infty) \quad  . $$

\paragraph{Axiom A0.} This axiom  states that for any $p\in \mathbb{N}$ and any $x \in X^{\omega}$ the dilatation 
$$\delta^{x}_{2^{p}} : U(x) \rightarrow V_{2^{p}}(x) $$
is a homeomorphism, the sets $U(x)$ and $\displaystyle V_{2^{p}}(x) $ are open and 
there is $A>1$ such that the ball centered in $x$ and radius $A$ is contained in $U(x)$. But this means 
that $\displaystyle U(x) = X^{\omega}$, because $\displaystyle X^{\omega} = B(x,1)$. 

Further, for any $p \in \mathbb{N}$ we have the inclusions: 
\begin{equation}
B(x, \frac{1}{2^{p}}) \subset \delta^{x}_{2^{p}} X^{\omega} \subset V_{2^{p}}(x) \quad . 
\label{a01}
\end{equation}

For any $\displaystyle p \in \mathbb{N}^{*}$ the associated dilatation  
$$\delta^{x}_{2^{-p}} : W_{2^{-p}}(x) \rightarrow B(x,B) = X^{\omega} \quad  , $$
is injective, invertible on the image. We suppose that $\displaystyle W_{2^{-p}}(x)$ is open, 
\begin{equation}
V_{2^{p}}(x) \subset W_{2^{-p}}(x) 
\label{a02}
\end{equation}
and that for all $\displaystyle p \in \mathbb{N}^{*}$ and $\displaystyle u \in X^{\omega}$ we have
$$\delta_{2^{-p}}^{x} \ \delta^{x}_{2^{p}} u \ = \ u \ .$$
We leave aside for the moment the interpretation of the technical condition before axiom A4.
 
\paragraph{Axioms A1 and A2.}  Nothing simplifies.

\paragraph{Axiom A3.} Because $d$ is an ultrametric distance and $\displaystyle X^{\omega}$ is compact, this axiom has very strong consequences, for a non degenerate dilatation structure. 

In this case the axiom A3 states that there is a non degenerate distance function  
$\displaystyle d^{x}$ on 
$\displaystyle X^{\omega}$ such that we have the limit 
\begin{equation}
\lim_{p \rightarrow \infty} 2^{p} d( \delta^{x}_{2^{p}} u, \delta^{x}_{2^{p}} v) = d^{x}(u,v) 
\label{a30}
\end{equation}
uniformly with respect to $\displaystyle x,u,v \in X^{\omega}$. 

We continue further with properties of weak dilatation structures. 

\begin{lema}
There exists $\displaystyle p_{0} \in \mathbb{N}$ such that for any $\displaystyle x, u, 
v \in X^{\omega}$ and for any $\displaystyle p \in \mathbb{N}$, $\displaystyle p \geq  p_{0}$, we have 
$$ 2^{p} d( \delta^{x}_{2^{p}} u, \delta^{x}_{2^{p}} v) = d^{x}(u,v) \quad  . $$
\label{l1}
\end{lema}

\paragraph{Proof.}
From the limit (\ref{a30}) and the non degeneracy of the distances $\displaystyle d^{x}$ we deduce that 
$$\lim_{p \rightarrow \infty} \log_{2} \left( 2^{p} d( \delta^{x}_{2^{P}} u, \delta^{x}_{2^{P}} v)\right)  =  
\log_{2} d^{x}(u,v) \quad , $$
uniformly with respect to $\displaystyle x,u,v \in X^{\omega}$, $u \not = v$. The right hand side term is finite and the sequence from the limit at the left hand side is included in $\mathbb{Z}$. Use this and 
the uniformity of the convergence to get the desired result. $\quad \square$
 
 In the sequel $\displaystyle p_{0}$ is the smallest natural number satisfying lemma \ref{l1}. 
 
 \begin{lema}
For any $\displaystyle x  \in X^{\omega}$ and for any $\displaystyle p \in \mathbb{N}$, $\displaystyle p \geq  p_{0}$, we have 
$$  \delta^{x}_{2^{p}} X^{\omega} = [x]_{p} X^{\omega}  \quad .  $$
Otherwise stated, for any $\displaystyle x, y \in X^{\omega}$, any $\displaystyle q  \in X^{*}$, 
$\displaystyle \mid q \mid \geq p_{0}$ there exists $\displaystyle w  \in X^{\omega}$ such that 
 $$  \delta^{qx}_{2^{\mid q \mid}} w = qy \quad , $$
 and for any $\displaystyle z \in X^{\omega}$ there is  $\displaystyle y \in X^{\omega}$ such that 
 $$  \delta^{qx}_{2^{\mid q \mid}} z = qy \quad . $$
 Moreover, for any $\displaystyle x  \in X^{\omega}$ and for any $\displaystyle p \in \mathbb{N}$, $\displaystyle p \geq  p_{0}$ the inclusions from (\ref{a01}), (\ref{a02}) are equalities. 
\label{l2}
\end{lema}

\paragraph{Proof.}
From the last inclusion in  (\ref{a01})  we get that for any $\displaystyle x, y \in X^{\omega}$, any $\displaystyle q  \in X^{*}$, 
$\displaystyle \mid q \mid \geq p_{0}$ there exists $\displaystyle w  \in X^{\omega}$ such that 
 $$  \delta^{qx}_{2^{\mid q \mid}} w = qy \quad . $$
 For the second part of the conclusion we use lemma \ref{l1} and axiom A1. From there we see that 
 for any $\displaystyle p \geq p_{0}$ we have 
 $$2^{p} d( \delta^{x}_{2^{p}} x, \delta^{x}_{2^{p}} u) =  2^{p} d( x, \delta^{x}_{2^{p}} u) = d^{x}(x,u) 
 \leq 1\quad  . $$
 Therefore $\displaystyle 2^{p} d( x, \delta^{x}_{2^{p}} u) \leq 1$, which is equivalent with the second part of the lemma. 
 
 Finally, the last part of the lemma has a similar proof, only that we have to use also the last part of axiom A0.  $\quad \square$

The technical condition before the axiom A4 turns out to be trivial. Indeed, from lemma \ref{l2} it follows that for any $\displaystyle p \geq p_{0}$, $p \in \mathbb{N}$, and any $\displaystyle x, u, v \in 
X^{\omega}$  we have $\displaystyle \delta^{x}_{2^{p}} u = [x]_{p} w$, $\displaystyle w \in X^{\omega}$.  It follows that  
$$\delta^{x}_{2^{p}} v \in [x]_{p} X^{\omega} = W_{2^{-p}} (x) = W_{2^{-p}} (  \delta^{x}_{2^{p}} u) \quad . $$

 \begin{lema}
For any $\displaystyle x, u, 
v \in X^{\omega}$  such that $\displaystyle 2^{p_{0}} d(x,u) \leq 1$,  $\displaystyle 2^{p_{0}} d(x,v) \leq 1$ we have 
$$d^{x}(u,v) = d(u,v) \quad . $$
Moreover, under the same hypothesis,  for any $\displaystyle p \in \mathbb{N}$ we have 
$$ 2^{p} d( \delta^{x}_{2^{p}} u, \delta^{x}_{2^{p}} v) = d(u,v) \quad  . $$
\label{l3}
\end{lema}

\paragraph{Proof.}
By lemma \ref{l1}, lemma \ref{l2} and axiom A2. Indeed, from lemma \ref{l1} and axiom A2,  for any $\displaystyle p \in \mathbb{N}$ and any 
$\displaystyle x, u', v' \in X^{\omega}$ we have 
$$d^{x}(u',v') = 2^{p_{0}+p} d( \delta^{x}_{2^{p+p_{0}}} u' ,  \delta^{x}_{2^{p+p_{0}}} v' ) = $$ 
$$ = 2^{p} \,   2^{p_{0}} d( \delta^{x}_{2^{p_{0}}}   \delta^{x}_{2^{p}}u', \delta^{x}_{2^{p_{0}}}   \delta^{x}_{2^{p}}v') = 2^{p} d^{x}( \delta^{x}_{2^{p}}u' , \delta^{x}_{2^{p}}v') \quad . $$
This is just the cone property for $\displaystyle d^{x}$. From here we deduce that for any 
$\displaystyle p \in \mathbb{Z}$ we have 
$$d^{x}(u', v') = 2^{p} d^{x}(\delta^{x}_{2^{p}}u' , \delta^{x}_{2^{p}}v') \quad . $$

If $\displaystyle 2^{p_{0}} d(x,u) \leq 1$,  $\displaystyle 2^{p_{0}} d(x,v) \leq 1$  then write $x = qx'$, $\displaystyle \mid q \mid = p_{0}$, and use lemma \ref{l2} to get the existence of $\displaystyle u', v' \in X^{\omega}$ such that 
$$\delta^{x}_{2^{p_{0}}} u' = u \quad , \quad \delta^{x}_{2^{p_{0}}} v' = v \quad . $$
Therefore, by lemma \ref{l1}, we have 
$$ d(u,v) = 2^{-p_{0}} d^{x}( u', v') = d^{x}( \delta^{x}_{2^{-p_{0}}} u' ,   \delta^{x}_{2^{-p_{0}}} v' ) = 
d^{x}(u,v) \quad . $$
The first part of the lemma is proven. For the proof of the second part write again 
$$2^{p} d( \delta^{x}_{2^{p}} u, \delta^{x}_{2^{p}} v) = 2^{p} d^{x}( \delta^{x}_{2^{p}} u, \delta^{x}_{2^{p}} v) = d^{x}(u,v) = d(u,v) \quad  \square $$

\subsection{Weak dilatation structures}

\begin{rk}
The space $\displaystyle X^{\omega}$ decomposes into a disjoint union of 
$\displaystyle 2^{p_{0}}$ balls which are isometric. There is no connection 
between the weak dilatation structures on these balls, therefore we shall 
study further only the case $\displaystyle p_{0} = 0$. 
\end{rk}

The purpose of this subsection is to find the  general form  of a weak dilatation structure   on $\displaystyle X^{\omega}$,  with  $\displaystyle p_{0} = 0$. 

\begin{defi}
A function $\displaystyle W: \mathbb{N}^{*} \times X^{\omega} \rightarrow 
Isom(X^{\omega})$  is smooth if for any $\varepsilon > 0$ there exists $\mu(\varepsilon) > 0$ such that 
for any $\displaystyle x, x' \in X^{\omega}$ such that $d(x,x')< \mu(\varepsilon)$ and for any $\displaystyle y \in X^{\omega}$  we have 
$$ \frac{1}{2^{k}} \, d( W^{x}_{k} (y) , W^{x'}_{k} (y) ) \leq \varepsilon \quad , $$
for an  $k$ such that  $\displaystyle d(x,x') <   1 / 2^{k} $.  
\label{defwsmooth}
\end{defi}

\begin{thm}
Let $\displaystyle (X^{\omega}, d, \delta)$ be a weak dilatation structure on $\displaystyle (X^{\omega}, d)$, where $d$ is the standard distance on $\displaystyle X^{\omega}$, such that $\displaystyle p_{0} = 0$. Then there exists a  smooth (according to definition \ref{defwsmooth}) function 
$$W: \mathbb{N}^{*} \times X^{\omega} \rightarrow Isom(X^{\omega}) \quad , \quad W(n,x) = W^{x}_{n} $$
such that  for any 
$\displaystyle q \in X^{*}$, $\alpha \in X$, $x, y \in X^{\omega}$ we have 
\begin{equation}
\delta_{2}^{q \alpha x} q \bar{\alpha} y = q \alpha \bar{x_{1}} W^{q \alpha x}_{\mid q \mid + 1} 
(y) \quad  . 
\label{eqtstruc}
\end{equation}
Conversely, to any smooth function  $\displaystyle W: \mathbb{N}^{*} \times X^{\omega} \rightarrow 
Isom(X^{\omega})$ is associated a weak dilatation structure  $\displaystyle (X^{\omega}, d, \delta)$, 
with $\displaystyle p_{0} = 0$, induced by  functions $\displaystyle \delta_{2}^{x}$, defined by 
$\displaystyle \delta_{2}^{x} x = x$ and otherwise by relation (\ref{eqtstruc}).
\label{tstruc}
\end{thm}
 
\paragraph{Proof.}
Let $\displaystyle (X^{\omega}, d, \delta)$ be a weak dilatation structure on $\displaystyle (X^{\omega}, d)$, such that $\displaystyle p_{0} = 0$. Any two different elements of $\displaystyle X^{\omega}$ can be written in the form $q \alpha x$ and $q \bar{\alpha} y$, with $\displaystyle q \in X^{*}$, $\alpha \in X$, $x, y \in X^{\omega}$. We also have 
$$d(q \alpha x , q \bar{\alpha} y) = 2^{- \mid q \mid} \quad . $$
From the following computation (using $\displaystyle p_{0} = 0$ and axiom A1): 
$$2^{-\mid q \mid -1} = \frac{1}{2} d(q \alpha x , q \bar{\alpha} y) = d( q \alpha x ,  \delta_{2}^{q \alpha x} q \bar{\alpha} y) \quad , $$
we find that there exists $\displaystyle w^{q \alpha x}_{\mid q \mid + 1}(y) \in X^{\omega}$ such that 
$$\delta_{2}^{q \alpha x} q \bar{\alpha} y = q \alpha w^{q \alpha x}_{\mid q \mid + 1}(y) \quad . $$
Further on, we compute: 
$$ \frac{1}{2} d(q \bar{\alpha} x , q \bar{\alpha} y) = d( \delta_{2}^{q \alpha x} q \bar{\alpha} x , \delta_{2}^{q \alpha x} q \bar{\alpha} y) = d( q \alpha w^{q \alpha x}_{\mid q \mid + 1}(x) , q \alpha w^{q \alpha x}_{\mid q \mid + 1}(y)) \quad . $$
From this equality we find that 
$$1 > \frac{1}{2} d(x,y) = d( w^{q \alpha x}_{\mid q \mid + 1}(x) , w^{q \alpha x}_{\mid q \mid + 1}(y)) \quad , $$
which means that the first letter of the word $\displaystyle w^{q \alpha x}_{\mid q \mid + 1}(y)$ 
does not depend on $y$, and is equal to the first letter of the word $\displaystyle w^{q \alpha x}_{\mid q \mid + 1}(x)$. Let us denote this letter by $\beta$ (which depends only on $q$, $\alpha$, $x$). Therefore we may write: 
$$w^{q \alpha x}_{\mid q \mid + 1}(y) = \beta W^{q \alpha x}_{\mid q \mid + 1}(y) \quad , $$
where the properties of the function $\displaystyle y \mapsto W^{q \alpha x}_{\mid q \mid + 1}(y)$ remain to be determined later. 

We go back to the first  computation in this proof: 
 $$2^{-\mid q \mid -1} = d( q \alpha x ,  \delta_{2}^{q \alpha x} q \bar{\alpha} y) = d( q \alpha x , q \alpha 
 \beta W^{q \alpha x}_{\mid q \mid + 1}(y)) \quad . $$
This shows that $\displaystyle \bar{\beta}$ is the first letter of the word $x$.  We proved the relation 
(\ref{eqtstruc}), excepting the fact that the function $\displaystyle y \mapsto W^{q \alpha x}_{\mid q \mid + 1}(y)$ is an isometry. But this is true. Indeed, for any $\displaystyle u,v \in X^{\omega}$ we have 
$$\frac{1}{2} d(q \bar{\alpha} u , q \bar{\alpha} v) = d( \delta_{2}^{q \alpha x} q \bar{\alpha} u , 
\delta_{2}^{q \alpha x} q \bar{\alpha} v) = d( q \alpha \bar{x_{1}} W^{q \alpha x}_{\mid q \mid + 1}(x) , q \alpha \bar{x_{1}} W^{q \alpha x}_{\mid q \mid + 1}(y)) \quad . $$
This proves the isometry property. 

The dilatations of coefficient $2$ induce all dilatations (by axiom A2). 
In order to satisfy the continuity 
assumptions from axiom A1, the function 
$\displaystyle W: \mathbb{N}^{*} \times X^{\omega} \rightarrow 
Isom(X^{\omega})$ has to be smooth in the sense of definition 
\ref{defwsmooth}.  Indeed, axiom A1 is equivalent to the fact that  
$\displaystyle \delta^{x'}_{2}(y')$ converges uniformly to $
\displaystyle \delta^{x}_{2}(y)$, as $d(x,x'), d(y,y')$ go to zero. 
There are two cases to study. 

Case 1: $d(x,x') \leq d(x,y)$, $d(y,y') \leq d(x,y)$. It means that $x = q \alpha q' \beta X$, $y = q \bar{\alpha} q" \gamma Y$, $x' = q \alpha q' \bar{\beta} X'$, $y' = q \bar{\alpha} q" \bar{\gamma} Y'$, 
with $d(x,y) = 1 / 2^{k}$, $k = \mid q \mid$ .

Suppose that $q' \not = \emptyset$. We compute then: 
$$ \delta_{2}^{x}(y) = q \alpha \bar{q'_{1}} W_{k+1}^{x}( q" \gamma Y) \quad , \quad \delta_{2}^{x'}(y') = 
q \alpha \bar{q'_{1}} W_{k+1}^{x'}(q" \bar{\gamma} Y') \quad . $$
All the functions denoted by a capitalized "W" are isometries, therefore we get the estimation: 
$$d(  \delta_{2}^{x}(y) ,  \delta_{2}^{x'}(y')) = \frac{1}{2^{k+2}} \, d(  W_{k+1}^{x}( q" \gamma Y) , W_{k+1}^{x'}(q" \bar{\gamma} Y')) \leq $$ 
$$ \leq \frac{1}{2^{k+2}} \, d( q" \gamma Y , q" \bar{\gamma} Y') + 
\, \frac{1}{2^{k+2}} \,  d( W_{k+1}^{x}(q" \gamma Y) , W_{k+1}^{x'}(q" \gamma Y)) = $$
$$ =  \frac{1}{2} d(y,y') + 
\, \frac{1}{2^{k+2}} \, d( W_{k+1}^{x}(q" \gamma Y) , W_{k+1}^{x'}(q" \gamma Y)) \quad . $$
We see that if $W$ is smooth in the sense of definition \ref{defwsmooth} then 
the structure $\delta$ satisfies the uniform continuity assumptions for this 
case. Conversely, if $\delta$ satisfies A1 then 
$W$ has to be smooth. 

If $q' = \emptyset$ then a similar computation leads to the same conclusion. 

Case 2: $d(x,x') > d(x,y) > d(y,y')$. It means that $x = q \alpha q' \beta X$, $x' = q \bar{\alpha} X'$, 
$y = q \alpha q' \bar{\beta} q" \bar{\gamma} Y$, $y' =  q \alpha q' \bar{\beta} q" \gamma Y'$, with 
$d(x,x') = 1 \ 2^{k}$, $k = \mid q \mid$ .  

We compute then: 
$$ \delta_{2}^{x}(y) = q \alpha q ' \beta \bar{X_{1}} W_{k+2+ \mid q' \mid}^{x}( q" \bar{\gamma} Y) \quad , \quad \delta_{2}^{x'}(y') = 
q \bar{\alpha} \bar{X'_{1}} W_{k+1}^{x'}(q' \bar{\beta} q" \gamma Y') \leq  $$
$$\leq \frac{1}{2^{k}} = d(x, x') \quad . $$
Therefore in his case the continuity is satisfied, without any supplementary constraints on the function $W$. 

The first part of the theorem is proven. 

For the proof of the second part of the theorem we start from the function 
$\displaystyle W: \mathbb{N}^{*} \times X^{\omega} \rightarrow Isom(X^{\omega})$. It is sufficient to prove  for any $\displaystyle 
x, y, z \in X^{\omega}$ the equality 
$$\frac{1}{2} d(y,z) = d(\delta_{2}^{x} y, \delta_{2}^{x} z) \quad .$$
Indeed, then we can construct the all dilatations from the dilatations of 
coefficient $2$ (thus we satisfy A2). All axioms, excepting A1, are satisfied. 
But A1 is equivalent with the smoothness of the function 
$W$, as we proved earlier.

Let us prove now the  before mentioned equality. 
If $y = z$ there is nothing to prove.  Suppose that $y \not = z$. The distance $d$ is ultrametric, therefore the proof splits in two cases. 

Case 1: $d(x,y) = d(x,z) > d(y,z)$. This is equivalent to $x = q \bar{\alpha} x'$, $y = q \alpha q' \beta y'$, 
$z = q \alpha q' \bar{\beta} z'$, with $\displaystyle q, q' \in X^{*}$, $\alpha, \beta \in X$, 
$\displaystyle x', y', z' \in X^{\omega}$. We compute: 
$$d( \delta_{2}^{x} y, \delta_{2}^{x} z) = d(\delta_{2}^{q \bar{\alpha} x'} q \alpha q' \beta y', \delta_{2}^{q \bar{\alpha} x'} q \alpha q' \bar{\beta} z') = $$
$$= d(q \bar{\alpha} \bar{x'_{1}} W^{x}_{\mid q \mid +1} ( q' \beta y') ,  q \bar{\alpha} \bar{x'_{1}} W^{x}_{\mid q \mid +1} ( q' \bar{\beta} z')) = 2^{-\mid q \mid - 1} d(  W^{x}_{\mid q \mid +1} ( q' \beta y'), 
 W^{x}_{\mid q \mid +1} ( q' \bar{\beta} z')) = $$
 $$ = 2^{-\mid q \mid - 1}  d(q' \beta y', q' \bar{\beta} z')) = \frac{1}{2} d(q \alpha q' \beta y' , q \alpha q' \bar{\beta} z') = \frac{1}{2} d(y,z) \quad  . $$

 Case 2: $d(x,y) = d(y,z) > d(x,z)$. If $x = z$ then we write $x = q \alpha u$, $y = q \bar{\alpha} v$ and we have 
 $$  d(\delta_{2}^{x} y, \delta_{2}^{x} z) = d( q \alpha \bar{u_{1}} W^{x}_{\mid q \mid +1} (v), q \alpha u) 
 = 2^{-\mid q \mid + 1} = $$
 $$ = \frac{1}{2} d(y,z) \quad . $$
If $x \not = z$ then we can write $z = q \bar{\alpha} z'$, $y = q \alpha q' \beta y'$, 
$x = q \alpha q' \bar{\beta} x'$, with $\displaystyle q, q' \in X^{*}$, $\alpha, \beta \in X$, 
$\displaystyle x', y', z' \in X^{\omega}$. We compute: 
$$d( \delta_{2}^{x} y, \delta_{2}^{x} z) = d(\delta_{2}^{q \alpha q' \bar{\beta} x'}  q \alpha q' \beta y' , 
\delta_{2}^{q \alpha q' \bar{\beta} x'} q \bar{\alpha} z') = $$
$$ = d(q \alpha q' \bar{\beta} \bar{x'_{1}} W^{x}_{\mid q \mid + \mid q' \mid +2}(y'), q \alpha \gamma 
W^{x}_{\mid q \mid + 1} (z')) \quad ,  $$
with $\gamma \in X$, $\displaystyle \bar{\gamma} = q'_{1}$ if $q' \not = \emptyset$, otherwise 
$\gamma = \beta$. In both situations we have 
$$d( \delta_{2}^{x} y, \delta_{2}^{x} z) =  2^{- \mid q \mid - 1} = \frac{1}{2} d(y,z) \quad . $$
The proof is done. $\quad \square$

\subsection{Self-similar  dilatation structures}

Let $\displaystyle (X^{\omega}, d, \delta)$ be a weak dilatation structure. There are induced 
dilatations structures on $\displaystyle 0X^{\omega}$ and  $\displaystyle 0X^{\omega}$. 

\begin{defi}
For any $\alpha \in X$ and $\displaystyle x, y \in X^{\omega}$ we define  $\displaystyle \delta_{2}^{\alpha, x} y$ by the relation 
$$\delta_{2}^{\alpha x} \alpha y = \alpha \, \delta_{2}^{\alpha, x} y \quad . $$
\end{defi}

The following proposition has a straightforward proof, therefore we skip it. 

\begin{prop}
If $\displaystyle (X^{\omega}, d, \delta)$ is a weak dilatation structure and $\alpha \in X$ then 
$\displaystyle (X^{\omega}, d, \delta^{\alpha})$ is a weak dilatation structure. 

If $\displaystyle (X^{\omega}, d, \delta')$ and $\displaystyle (X^{\omega}, d, \delta")$ are weak dilatation structures then $\displaystyle (X^{\omega}, d, \delta)$ is a weak dilatation structure, where 
$\delta$ is uniquely defined by $\displaystyle \delta^{0} = \delta'$, $\displaystyle \delta^{1} = \delta"$. 
\end{prop}

The previous proposition justifies the next definition. 

\begin{defi}
 A weak dilatation structure $\displaystyle (X^{\omega}, d, \delta)$ is self-similar if for any $\alpha \in X$ and $\displaystyle x,y \in X^{\omega}$ we have 
 $$\delta_{2}^{\alpha x} \, \alpha y = \alpha \, \delta_{2}^{x} y \quad . $$
\end{defi}

\begin{prop}
Let $\displaystyle (X^{\omega}, d, \delta)$ be a self-similar weak dilatation structure and 
$\displaystyle W: \mathbb{N}^{*} \times X^{\omega} \rightarrow Isom(X^{\omega})$ the function 
associated to it, according to theorem \ref{tstruc}.  Then there exists a function $W: X^{\omega} \rightarrow Isom(X^{\omega})$ such that: 
\begin{enumerate}
\item[(a)] for any $\displaystyle q \in X^{*}$ and any $\displaystyle x \in X^{\omega}$ we have 
$$W_{\mid q \mid + 1}^{q x} = W^{x} \quad , $$
\item[(b)] there exists $C>0$ such that for any $\displaystyle x, x', y \in X^{\omega}$ and for any 
$\lambda > 0$, if $d(x,x') \leq \lambda$ then 
$$d(W^{x}(y) , W^{x'}(y)) \leq C \lambda \quad . $$
\end{enumerate}
\end{prop}

\paragraph{Proof.}
We define $\displaystyle W^{x} = W^{x}_{1}$ for any $\displaystyle x \in X^{\omega}$ We want to prove 
that this function satisfies (a), (b). 

(a) Let $\displaystyle \beta \in X$ and any $\displaystyle x, y \in X^{\omega}$, $x = q \alpha u$, $y = q \bar{\alpha} v$. By self-similarity we obtain: 
$$\beta q \alpha \bar{u_{1}} W_{\mid q \mid + 2}^{\beta x} (v) = \delta_{2}^{\beta x} \beta y =  \beta \delta_{2}^{x} y =  \beta q \alpha \bar{u_{1}} W_{\mid q \mid + 1}^{x} (v) \quad . $$
We proved that 
$$ W_{\mid q \mid + 2}^{\beta x} (v) = W_{\mid q \mid + 1}^{x} (v)  $$
for any $\displaystyle x, v \in X^{\omega}$ and $\beta \in X$ This implies (a).  

(b) This is a consequence of smoothness,  in the sense of definition \ref{defwsmooth}, of the function 
$\displaystyle W: \mathbb{N}^{*} \times X^{\omega} \rightarrow Isom(X^{\omega})$.  Indeed, 
$\displaystyle (X^{\omega}, d, \delta)$ is a  weak dilatation structure, therefore by theorem \ref{tstruc} 
the  previous mentioned function is smooth. 

By (a) the smoothness condition becomes: for any $\varepsilon > 0$ there is $\mu(\varepsilon) > 0$ such that for any $\displaystyle y \in X^{\omega}$, any $k \in \mathbb{N}$ and  any $x, x' \in X^{\omega}$,  
if $d(x,x') \leq 2^{k} \mu(\varepsilon)$ then 
$$d( W^{x}(y), W^{x'} (y)) \leq 2^{k} \varepsilon \quad . $$
Define then the modulus  of continuity: for any $\varepsilon > 0$ let $\bar{\mu}(\varepsilon)$ be given by 
$$\bar{\mu}(\varepsilon) = \sup \left\{ \mu \, \mbox{ : } \forall x, x', y \in X^{\omega} \, \, d(x,x')\leq \mu \Longrightarrow d( W^{x}(y), W^{x'} (y)) \leq \varepsilon \right\} \quad .$$
We see that the modulus of continuity $\bar{\mu}$ has the property 
$$\bar{\mu}(2^{k} \varepsilon)  = 2^{k} \bar{\mu}(\varepsilon) $$
for any $k \in \mathbb{N}$. Therefore there exists $C> 0$ such that 
$\displaystyle \bar{\mu}( \varepsilon) = C^{-1} \varepsilon$ for any $\displaystyle \varepsilon = 1/2^{p}$, $p \in \mathbb{N}$. The point (b) follows immediately. \quad $\square$

\section{The Cantor set again}

\begin{thm}
Let $\displaystyle (X^{\omega}, d, \delta)$ and $\displaystyle (X^{\omega}, d, \bar{\delta})$ be two linear 
and self-similar weak dilatation structures, with $\displaystyle p_{0} = 0$. Suppose that 
 there are two different points $\displaystyle x_{0}, x_{1} \in X^{\omega}$ such that 
$\displaystyle \delta_{2}^{x_{i}} = \bar{\delta}_{2}^{x_{i}}$, $i=0,1$.  Then  a dense 
set $\displaystyle B \subset X^{\omega}$ exists such that for any $x \in B$ there is 
a positive radius $r(x)>0$ such that for any $\displaystyle y \in X^{\omega}$ with 
$d(x,y)\leq r(x)$ we have: 
$$\delta_{2}^{x} y \ = \ \bar{\delta}^{x}_{2} y \quad . $$ 
\label{tunic}
\end{thm}

\paragraph{Proof.} 
We use the self-similarity hypothesis to restrict to the case of  
$\displaystyle x_{0},  x_{1} \in X^{\omega}$, such that  $\displaystyle d(x_{0}, x_{1}) = 1$.

Let $q\in X^{*}$, $q \not = \emptyset$, a non empty word. To $q= q_{1} ... q_{n}$ is associated 
$$f(q) = \delta_{2}^{x_{q}} x_{0} = \delta_{2}^{x_{q_{1}}} \, \delta_{2}^{x_{q_{2}}} \, ... \delta_{2}^{x_{q_{n}}}  
(x_{0}) \, \in X^{\omega} \quad . $$
 
 We shall use the next lemma. 
 
\begin{lema}
Let $\displaystyle F = \left\{ 0 \right\} \cup  X^{*}1$ be the set of all non empty finite words  which 
either are equal to $0$  or they end with 1. Then the restriction of the previously defined function $f$ to  $$\displaystyle f: F \rightarrow X^{\omega}$$ is injective and the image of $f$ is dense in $\displaystyle X^{\omega}$.  Moreover, 
$\displaystyle f(X^{*}\setminus \left\{ \emptyset \right\}) = f(F)$. 
\label{le1}
\end{lema}

The dilatation structure  $\displaystyle (X^{\omega}, d, \delta)$ is linear by hypothesis, which 
implies that for  any $x = f(q)$, $q\in F$, $\mid q \mid = n$, and 
$$\displaystyle u \in \delta_{2}^{x_{q}} X^{\omega} =  [f(q)]_{n} X^{\omega}$$ 
the value of the dilatation $\displaystyle \delta_{2}^{f(q)} (u)$ is uniquely determined by the relation: 
$$\delta_{2}^{f(q)} (u)  = \delta_{2}^{x_{q}} \, \delta_{2}^{x_{0}} \, \left( \delta_{2}^{x_{q}} \right)^{-1} (u) \quad . $$
Form the hypothesis we deduce that  in particular for $q, u$ as previously  we have 
$$\delta_{2}^{f(q)} (u) = \bar{\delta}_{2}^{f(q)} (u) \quad  . $$
From the self-similarity of the dilatation structures we deduce that we have the same equality for all 
$u \in X^{\omega}$ and for all $q$. From the continuity of dilatation structures and the density result 
stated in lemma \ref{le1} we get the result. \quad $\square$

\paragraph{Proof of Lemma \ref{le1}.}
We shall suppose that $\displaystyle x_{0} = (0) = 0000...$ and 
$\displaystyle x_{1}= (1) = 1111...$, only for expository reasons.  
The proof which follows may be easily (but with longer notations) 
adapted to the general case. 

Let $\displaystyle q_{1}, q_{2} \in F$ such that  
$\displaystyle q_{1}= z \alpha q_{1}', \, 
q_{2}= z \bar{\alpha} q_{2}'$. Then 
$$d(f(q_{1}), f(q_{2})) = d( \delta_{2}^{x_{z}} \delta_{2}^{x_{\alpha}} \delta_{2}^{x_{q_{1}'}} (x_{0}) , 
 \delta_{2}^{x_{z}} \delta_{2}^{x_{\bar{\alpha}}} \delta_{2}^{x_{q_{2}'}} (x_{0})) = $$ 
 $$ = \frac{1}{2^{\mid z \mid}} d ( \delta_{2}^{x_{\alpha}} \delta_{2}^{x_{q_{1}'}} (x_{0}) ,  \delta_{2}^{x_{\bar{\alpha}}} \delta_{2}^{x_{q_{2}'}} (x_{0}))  = \frac{1}{2^{\mid z \mid}} d ( \delta_{2}^{x_{0}} (u) ,  \delta_{2}^{x_{1}} (v)) \quad , $$
 for certain $\displaystyle u,v \in  X^{\omega}$. It follows that 
 $$d(f(q_{1}), f(q_{2})) = \frac{1}{2^{\mid z \mid}} \quad . $$
 Suppose now that $\displaystyle q_{1}, q_{2} \in F$ , $\displaystyle q_{1} = q_{2} q$, with $q\not = \emptyset$. Then $q \in F$, $q \not = 0$  and 
 $$ d(f(q_{1}), f(q_{2})) = d( \delta_{2}^{x_{q_{1}}}  \delta_{2}^{x_{q}} (x_{0}) , 
  \delta_{2}^{x_{q_{1}}} (x_{0})) = \frac{1}{2^{\mid q_{1} \mid}} d( \delta_{2}^{x_{q}} (x_{0}) , x_{0}) \quad . 
  $$
  We want to prove that for any $q\in F$ $\displaystyle  d( \delta_{2}^{x_{q}} (x_{0}) , x_{0}) \not = 0$. 
  We know that  $q= q' 1$ and we use a ping-pong type reasoning. Notice that $f(1) \in 10X^{\omega}$. 
  We shall prove that for any $\displaystyle q \in X^{*}$, $q \not = 0$, we cannot have 
  $$(0) = \delta_{2}^{x_{q}} f(1) \quad . $$  
  Indeed, for any $\displaystyle p, r\in \mathbb{N}^{*}$ we have 
  $$\delta_{2^{p}}^{(0)} 0^{r}1X^{\omega} \subset  0^{p+r}1 X^{\omega} \quad . $$
Remark  that if we apply $\displaystyle \delta_{2}^{(1)}$ to $\displaystyle X^{\omega}$ then we get 
  $\displaystyle 1X^{\omega}$. Therefore after application of a finite string of $\displaystyle \delta_{2}^{(0)},  \delta_{2}^{(1)}$ to a word starting with $10$, we shall always get a word containing $1$. This proves the first part of the lemma. 
  
  For the second part remark that 
  $$X^{\omega} = \delta_{2}^{(0)} X^{\omega} \, \cup \, \delta_{2}^{(1)} X^{\omega} \quad $$
and use the Hutchinson theorem \ref{thutch}. The lemma is proven. \quad $\square$

The standard IFS of the Cantor set $\displaystyle X^{\omega}$ is given by the pair 
of transformations $\displaystyle \phi_{\alpha}(x) = \alpha x$, $\alpha = 0,1$.  
These two transformations are contractions of coefficient 1/2, therefore we may think 
about them as dilatations of coefficient $2$ based at $(0) = 0000...$ and $(1)=1111...$ respectively. 

Suppose that   $\displaystyle (X^{\omega}, d, \delta)$ is a linear, self-similar 
dilatation structure with the property that for any $\displaystyle x \in X^{\omega}$: 
\begin{equation}
 \delta_{2}^{(0)} x = \phi_{0}(x) = 0x \quad , 
\quad  \delta_{2}^{(1)} x = \phi_{1}(x) = 1x \quad .
\label{eqned1}
\end{equation}
Then, according to theorem \ref{tunic}, the dilatation structure is uniquely determined 
on a dense subset of $\displaystyle X^{\omega}$, in the sense explained in the conclusion 
of the theorem. 

Is this implying that there is only one such dilatation structure? 
(According to examples, there is at least one such dilatation structure, coming from 
dyadic integers). 

We shall even put more demand on the dilatation structure 
$\displaystyle (X^{\omega}, d, \delta)$ by asking further that for any $\displaystyle 
n \in \mathbb{N}^{*}$, any  
finite word of length $n$  $\displaystyle q \in X^{n}$,  
the composition of functions $\displaystyle \phi_{q_{n}} ... \phi_{q_{1}}$ 
is a dilatation of coefficient $\displaystyle 2^{n}$, based at 
$(q) = qqqq...$, that is:
\begin{equation} 
\phi_{q_{n}} ... \phi_{q_{1}} = \delta^{(q)}_{2^{n}} \quad .
\label{eqned2}
\end{equation}
Is there more than one linear, self similar dilatation structure which satisfies 
relations (\ref{eqned1}), (\ref{eqned2})?

\begin{prop}
The standard IFS of the Cantor set does not uniquely determine a 
linear and self-similar dilatation structure satisfying  (\ref{eqned1}), (\ref{eqned2}).
There exist two different and non equivalent strong linear dilatation structures which 
satisfy these two relations.
\label{corunic}
\end{prop}

\paragraph{Proof.}
There is a strong, linear and self-similar dilatation structure which has these two 
contractions $\displaystyle \phi_{\alpha}$, $\alpha=0,1$, as dilatations. 
This is coming from the 
interpretation of $\displaystyle X^{\omega}$ as the set of dyadic integers. 
Under this identification 
we have 
$$\delta_{2}^{x} y = x + 2(y-x) = 2 y - x \quad . $$
The dilatation structure defined 
like this is strong, linear and self-similar. Let us check this. 

Start with linearity: for any $x, u, v \in X^{\omega} \equiv \mathbb{Z}_{2}$ the following is true
$$\delta_{2}^{x} \delta_{2}^{u} (v) = \delta_{2}^{\delta_{2}^{x} u} \delta_{2}^{x} v \quad . $$
Indeed, the left hand side equals: 
$$\delta_{2}^{x} \delta_{2}^{u} (v) = \delta_{2}^{x} (2v - u) = 2(2v-u) - x \quad , $$
and the right hand side is equal to 
$$\delta_{2}^{\delta_{2}^{x} u} \delta_{2}^{x} v = \delta_{2}^{2u - x} (2v - x) = 2(2v -x) - (2u -x) = 2(2v-u) - x \quad . $$

Concerning the self-similarity: let $\displaystyle q\in X^{*}$, non empty finite word of lenght $n$, and 
$\displaystyle u,v \in X^{\omega}$. Then we may identify $q$ with the element of $\displaystyle 
\mathbb{Z}_{2}$ described by the word $q(0)$, and we have 
$$qu = q + 2^{n}u \quad , \quad qv = q + 2^{n} v \quad , $$
therefore the self-similarity comes from the string of equalities: 
$$\delta_{2}^{qu} qv = 2(q + 2^{n}v) - (q + 2^{n}u) = q + 2^{n} ( 2v - u) = q \delta_{2}^{u} v \quad . $$
The dilatation structure is strong because we have here a conical group. 

We check now the relations (\ref{eqned1}), (\ref{eqned2}). These are equivalent with: 
 for any $\displaystyle x \in 
X^{\omega}$ and for any  $\displaystyle 
n \in \mathbb{N}^{*}$, any  
finite word of length $n$  $\displaystyle q \in X^{n}$, we have: 
$$\delta^{(q)}_{2^{n}} x = qx \quad .$$
This is true for our dilatation structure. Indeed, when we identify $\displaystyle X^{\omega}$ 
with the dyadic integers, the word $(q)$ is identified with 
$$(q) = \sum_{i=1}^{n} q_{i} 2^{i-1}\left( \sum_{k=0}^{\infty} 2^{k} \right) = 
\frac{q}{1-2^{n}} \quad , $$\
(here we used the identification $\displaystyle q = \sum_{i=1}^{n} q_{i} 2^{i-1}$).  Therefore 
$$ \delta^{(q)}_{2^{n}} x = \frac{q}{1-2^{n}} + 2^{n} \left( x - \frac{q}{1-2^{n}} \right) = 
 q + 2^{n} x = qx \quad . $$

The second strong linear, self-similar dilatation structure which satisfies  the 
relations (\ref{eqned1}), (\ref{eqned2}) is the 
dilatation structure associated to the infinite dihedral group $\displaystyle 
\mathbb{D}_{\infty}$. This is the group $\displaystyle (X^{\omega},*)$ 
with the following group operation: for any $\displaystyle x,y \in X^{\omega}$ the sum 
$\displaystyle x*y \in X^{\omega}$ is the word with the letters 
$\displaystyle (x*y)_{i} = x_{i} + y_{i}  \mod  2$. Remark that for 
any $\displaystyle x \in X^{\omega}$ we have $\displaystyle x^{-1} = x$. 

Consider the group morphism: $\displaystyle \delta_{2} x = 0x$. This induces the dilatation structure 
$$ \bar{\delta}_{2}^{x} (y) = x * \delta_{2}(x*y) \quad . $$
This is a strong linear dilatation structure. It is also self-similar and it satisfies 
(\ref{eqned1}), (\ref{eqned2}) by straightforward computations.

These two dilatation structures are not equivalent. Indeed, both dilatation structures 
come from conical groups. If they are equivalent then the groups (dyadic integers 
with addition and the infinite dihedral group) have to be isomorphic. But these two groups are 
not isomorphic therefore the dilatation structures cannot be equivalent.  \quad $\square$

In contrast, on $\displaystyle \mathbb{R}$ with euclidean distance, up to 
equivalence, there is only one strong linear  dilatation structure coming from a 
conical group without small subgroups.

\end{document}